\documentclass[english, aps, showpacs, preprint, nobibnotes]{revtex4-1}

\usepackage{amsmath,amssymb}

\usepackage[utf8x]{inputenc}

\usepackage[T1]{fontenc}
\usepackage{float}
\usepackage{mathtools}
\usepackage{amsmath}
\usepackage{amssymb}
\usepackage{graphicx}
\usepackage{nicefrac}
\usepackage{resizegather}
\usepackage{breakcites}

\begin{document}
	
	\title{A Center Manifold Reduction Technique \\for a System of Randomly Coupled Oscillators}
	\date{\today}
	
	\author{Dimitrios Moirogiannis}
	\email{dmoirogi@gmail.com}
	\author{Keith Hayton}
	\email{khayt86@gmail.com }
	\author{Marcelo Magnasco}
	\affiliation{Center for Studies in Physics and Biology, The Rockefeller University}

	\begin{abstract}
		
	In dynamical systems theory, a fixed point of the dynamics is called \textit{nonhyperbolic} if the linearization of the system around the fixed point has at least one eigenvalue with zero real part. The \textit{center manifold existence theorem} guarantees the local existence of an invariant subspace of the dynamics, known as a center manifold, around such nonhyperbolic fixed points. A growing number of theoretical and experimental studies suggest that some neural systems utilize nonhyperbolic fixed points and corresponding center manifolds to display complex, nonlinear dynamics and to flexibly adapt to wide-ranging sensory input parameters. In this paper, we present a technique to study the statistical properties of high-dimensional, nonhyperbolic dynamical systems with random connectivity and examine which statistical properties determine both the shape of the center manifold and the corresponding reduced dynamics on it. This technique also gives us constraints on the family of center manifold models that could arise from a large-scale random network. We demonstrate this approach on an example network of randomly coupled damped oscillators.      
		
	\end{abstract}

	\maketitle

	\section{Introduction}	
	The \textit{manifold hypothesis} proposes that phenomena in the physical world often lie close to a lower dimensional manifold embedded in the full high-dimensional space \cite{geron2017}. As technological advancements have allowed the collection of large and often complex data sets, the manifold hypothesis has spurred the creation of a number of \textit{manifold learning}, or \textit{manifold reduction}, techniques \cite{belkin2003,carlsson2009,genovese2012,dasgupta2008,hastie1989,kambhatla1994,kegl2000,narayanan2009,niyogi2008,perrault2012,roweis2000,smola2001,tenenbaum2000,weinberger2006}. In this era of ``Big Data'', dimensionality reduction techniques are often part of the initial step in the data processing pipeline. Reducing dimensions can speed up model fitting, aid in visualization, and avoid the pitfalls of overfitting in high-dimensional data sets \cite{geron2017}.
	
	In the field of neuroscience, where large collections of high-dimensional data are often the focus of research, a growing number of studies suggest that neural systems lie close to a dynamical systems structure known as a \textit{center manifold}. These studies include entire hemisphere ECoG recordings \cite{Solovey2012, Alonso2014, solovey2015loss}, experimental studies in premotor and motor cortex \cite{churchland2012}, theoretical \cite{seung1998continuous} and experimental studies \cite{seung2000stability} of \textit{slow manifolds} (a specific case of center manifolds) in oculomotor control, slow manifolds in decision making \cite{machens2005}, Hopf bifurcation \cite{poincare1893, Hopf1942, Andronov2013} in the olfactory system \cite{freeman2005metastability} and cochlea \cite{choe1998model, eguiluz2000essential, camalet2000auditory, kern2003essential, Duke2003, Magnasco2003, Hayton2018b}, a nonhyperbolic model of primary visual cortex \cite{Hayton2018a}, and theoretical work on regulated criticality \cite{bienenstock1998regulated}. 
	
	To understand center manifolds, we now review some basic ideas and definitions from dynamical systems theory. First, the classical approach to studying behavior in the neighborhood of an equilibrium point is to examine the eigenvalues of the Jacobian of the system at this point. If all eigenvalues have nonzero real part, then the equilibrium point is called \textit{hyperbolic}. In this case, according to the \textit{Hartman-Grobman theorem}, the dynamics around the point is topologically conjugate to the linearized system determined by the Jacobian \cite{Grobman1959,Hartman1960a,Hartman1960b}. This implies that solutions of the system near the hyperbolic fixed point exponentially decay (or grow) with time constants determined by the real part of the eigenvalues of the linearization. Thus, locally, the nonlinearities in the system do not play an essential role in determining the dynamics.
	
	This behavior around hyperbolic points is in stark contrast to the dynamics around \textit{nonhyperbolic} fixed points, where there exists at least one eigenvalue of the Jacobian with zero real part \cite{Izhikevich2007,hoppensteadt2012weakly,wiggins2003introduction}. The linear space spanned by the eigenvectors of the Jacobian corresponding to the eigenvalues on the imaginary axis (\textit{critical modes}) is called the \textit{center subspace}. In the case of nonhyperbolicity, the Hartman-Grobman theorem does not apply, the dynamics are not enslaved by the exponent of the Jacobian, and nonlinearities play a crucial role in determining dynamical properties around the fixed point. The classical approach to studying this nonlinear behavior around nonhyperbolic fixed points is to investigate the reduced dynamics on an invariant subspace called a \textit{center manifold}. The \textit{center manifold existence theorem} \cite{Kelley1967,Carr2012} guarantees the local existence of this invariant subspace. It is tangent to the \textit{center subspace}, and its dimension is equal to the number of eigenvalues on the imaginary axis. If it is the case that all eigenvalues of the linearization at a nonhyperbolic fixed point have zero or negative real part (no unstable modes), as is often the case in physical systems, then for initial conditions near the fixed point, trajectories exponentially approach solutions on the center manifold \cite{Carr2012}, and instead of studying the full system, we can study the reduced dynamics on the center manifold.
	
	Dynamics on center manifolds around nonhyperbolic fixed points are complex and can give rise to interesting nonlinear features \cite{Adelmeyer1999, wiggins2003introduction, Kuznetsov2013}. In general, the greater the number of eigenvalues on the imaginary axis, the more complex the dynamics could be. First, since the dynamics are not enslaved by the exponent of the Jacobian, nonlinearities and input parameters play a crucial role in determining dynamical properties such as relaxation timescales and correlations \cite{yan2012input, Hayton2018a}. Moreover, activity perturbations on the center manifold neither damp out exponentially nor explode, but instead, perturbations depend algebraically on time. If the center manifold is also delocalized, which can occur even with highly local connections between individual units (units not connected on the original network can be connected on the reduced system on a center manifold), activity perturbations can then also propagate over large distances. This is in stark contrast to the behavior in stable systems where perturbations are damped out exponentially and information transfer can only occur on a timescale shorter than the timescale set by the system's exponential damping constants. 
	
	Nonhyperbolic equilibrium points are also referred to as \textit{critical points} \cite{ ShnoScholarpedia, Izhikevich scholarpedia} and the resulting dynamics as \textit{dynamical criticality}. Dynamical criticality is distinct from \textit{statistical criticality} \cite{beggs2012being}, which is related to the statistical mechanics of second-order phase transitions. It has been proposed that neural systems \cite{chialvo2010emergent}, and more generally biological systems \cite{mora2011biological}, are statistically critical in the sense that they are poised near the critical point of a phase transitions \cite{da1998criticality,fraiman2009ising}. Statistical criticality is characterized by power law behavior such as avalanches \cite{beggs2003neuronal,levina2007dynamical,gireesh2008neuronal} and long-range spatiotemporal correlations \cite{eguiluz2005scale,kitzbichler2009broadband}. While both dynamical criticality and statistical criticality have had success in neuroscience, their relation is still far from clear \cite{magnasco2009self,mora2011biological,kanders2017avalanche}.
	
	In light of the growing evidence that center manifolds play a crucial role in neural dynamics, we feel that it is prudent to develop dimensionality reduction techniques aimed at reducing dynamics which contain a prominent center subspace to more simple low-dimensional dynamics on center manifolds. In this paper, we explore center manifold reduction as a statistical nonlinear dimensionality reduction technique on an example network having random connectivity as its base. The random connectivity is inspired by the Sompolinsky family of neural network models \cite{Sompolinsky1998chaos}, which assume an interaction structure given by a Gaussian random matrix scaled to have a spectrum in the unit disk. This family of neural networks has led to a number of useful results, including a flurry of recent publications \cite{Stern2014Dynamics, Lalazar2016Tuning, Rajan2016Recurrent, Ostojic2014two, Kadmon2015transition, Landau2016Impact, Sompolinsky2014Computational, Engelken2016reanalysis, Harish2015Asynchronous}. 
	
	We study which statistical properties of the random connectivity are essential in determining both the shape of the center manifold and the corresponding reduced dynamics on it. We have previously \cite{Moirogiannis2017} presented a naive linear approach of this calculation in another paper, but here we will include the nonlinearities of the center manifolds as well. Using the center manifold reduction algorithm, we find that the resulting equations depend crucially on the structure and the statistics of the eigenvectors of the connectivity matrix. Even though the statistical distribution of eigenvalues of a broad class of random matrices has been extensively characterized \cite{Tao2008circular}, the statistical distribution of the eigenvectors is not well understood \cite{Chalker1998Eigenvector, Tao2012Random, O'Rourke2016Eigenvectors}. Our analysis shows that the reduced dynamics on the center manifold, and thus the collective dynamics of the full system, can be drastically different from the dynamics of the individual subunits that are randomly coupled. We also find interesting constraints on the family of center manifold models that could arise from a large-scale random network. 
	
\begin{figure}[bth]
			\begin{centering}
				\includegraphics[scale=0.6]{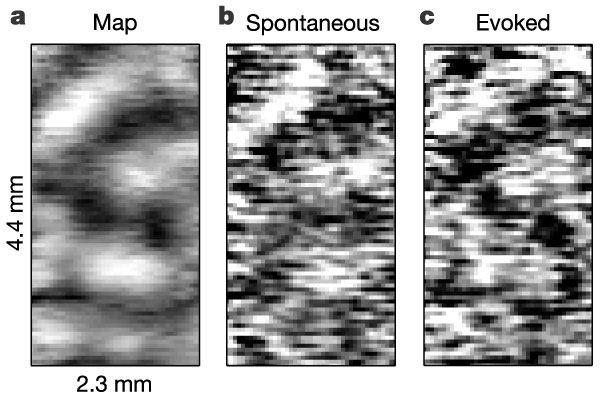}
				\par\end{centering}
			
			\begin{centering}
				\caption{\label{fig:spontaneous} \textbf{Reprinted by permission from Springer Nature: Nature. Spontaneously emerging cortical representations of visual attributes. Kenet \textit{et al}. COPYRIGHT 2003.} Activity corresponding closely to orientation maps spontaneously arises in the absence of input in an area with orientation selectivity. a) An orientation map of vertical orientation from cat area 18, where most cells are selective for stimulus orientation, obtained by using voltage sensitive dye imaging. b) A map obtained in a single frame from a spontaneous recording session. c) A single frame from an evoked session using the same orientation as for the map.}
				\par\end{centering}
			
		\end{figure}

	In our example, we only consider spontaneous activity; the network is not driven by an external forcing. A number of studies have highlighted the importance of spontaneous activity in cortex \cite{arieli1996,petersen2003,Fox2005,Fox2006}. For example, the spontaneous activity of individual units in primary visual cortex is strongly coupled to global patterns evincing the underlying functional architecture \cite{Tsodyks1999, Kenet2003}: during spontaneous activity in the absence of input, global modes are transiently activated that are similar to modes engaged when inputs are presented (Fig. \ref{fig:spontaneous}). We therefore feel that it is important to understand center manifold reduction in the context of spontaneous activity, as we do in this paper.

	\begin{figure}[bth]
		\begin{centering}
			\includegraphics[scale=0.6]{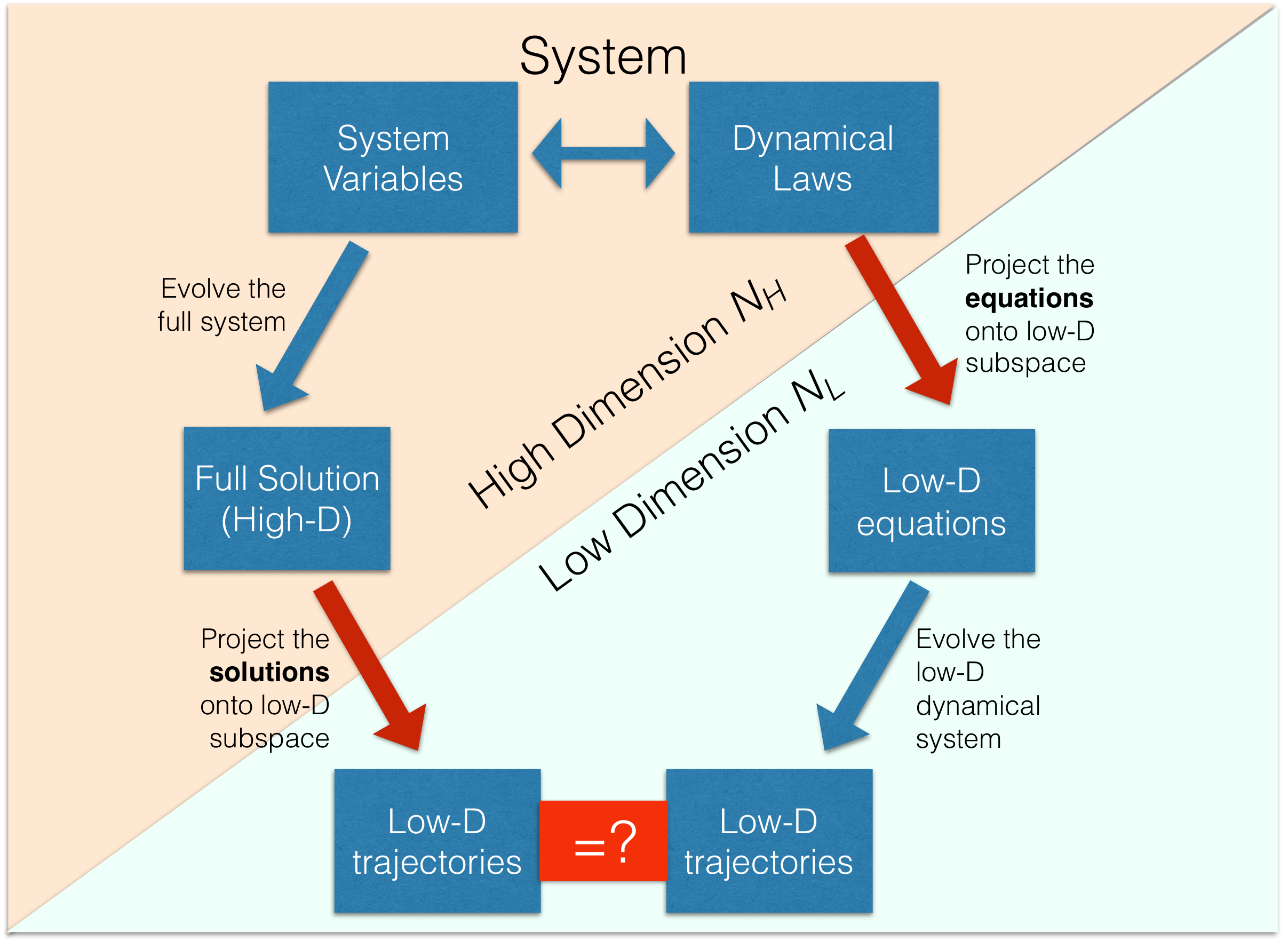}
			\par\end{centering}
		
		\begin{centering}
			\caption{\label{fig:schematics} \textbf{Schematics of Dimensionality Reduction}}
			\par\end{centering}
		
	\end{figure}

	The general schematics of our dimensionality reduction approach is illustrated in (Fig. \ref{fig:schematics}). At the top left we have the original system variables, which lives in dimension $N_H$, ruled by a high-dimensional dynamical law. At the bottom right, we have a low-dimensional space, the reduced variables, which live in dimension $N_L$. On the left branch, we first perform a detailed, costly and precise simulation of every variable in the system, and  then we project this high-dimensional solution down into the lower-dimensional space. On the right branch, instead, we first project down the original high-dimensional vector field to obtain the reduced dynamics on the lower-dimensional subspace. We then evolve the reduced dynamics in its low-dimensional space. If both branches give rise to approximately the same low-D trajectories, then the reduction is well defined.

	The structure of the paper is as follows. We first introduce some general definitions and the network structure of the specific example we will be using. We then review how the linear projection of the dynamics on the center space can be used for as a reduction approximation for certain simple cases \cite{Moirogiannis2017}, but also how it can fail for a wider range of parameters of the system. Next, we introduce the main focus of the paper, the proposed center manifold reduction technique, which can be thought as a nonlinear generalization of the linear projection on the center subspace. Finally, using the resulting approximations as numerical Ansatz to correct both the reduced dynamics equations in the case where the naive approximation fails as well as to describe the dynamics of individual units using the reduced dynamics.
	
	 \section{Definitions}
	
	Let us now establish the notation to be used in further detail, shown schematically in Fig. \ref{fig:schematics2}:

		\begin{figure}[bth]
		\begin{centering}
			\includegraphics[scale=0.6]{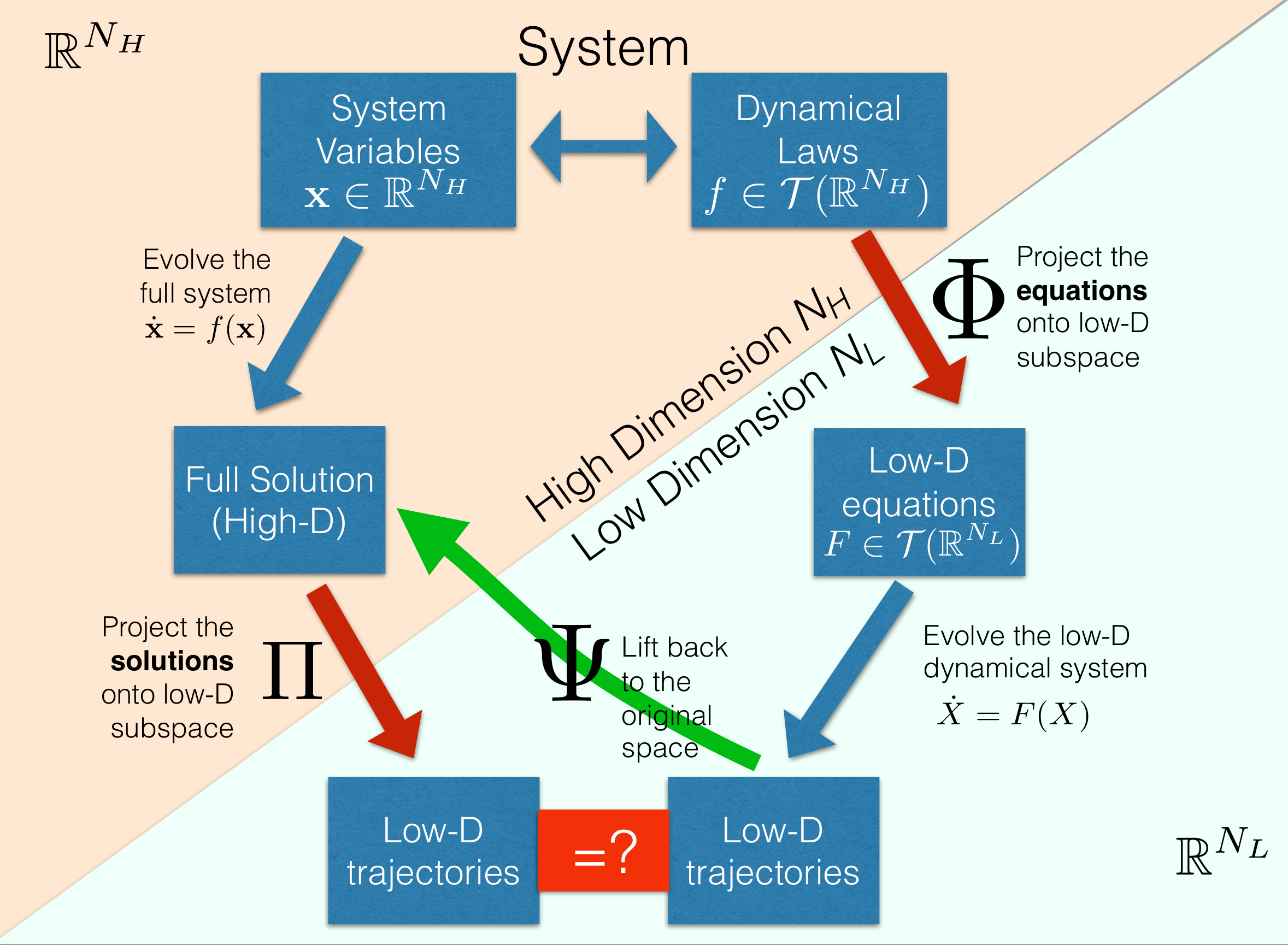}
			\par\end{centering}
		
		\begin{centering}
			\caption{\label{fig:schematics2} \textbf{Operator Definitions}}
			\par\end{centering}
		
	\end{figure}

	Again, at the top left we have our original full system: many variables $ x \in {\mathbb R}^{N_H}$ ruled by a high-dimensional dynamical law $\dot  x = f( x)$ given by a vector field $f$ in $ {\mathbb R}^{N_H}$, stated as $f \in {\mathcal T}({\mathbb R}^{N_H})$. At the bottom right we have a low-dimensional space, the reduced variables, which live in ${\mathbb R}^{N_L}$. 
	On the left branch we first evolve in $N_H$ and  then we use a projection operator  $\Pi$ to project the high-dimensional solution down into the lower-dimensional space  On the right branch, instead, we  first project down the high-dimensional vector field to obtain a reduced vector field $F$ on the lower-dimensional subspace: $F \in {\mathcal T}({\mathbb R}^{N_L})$. If we call the projector operator for the equations $\Phi$; then $F = \Phi\{f\}$. Then we evolve this simpler dynamics $\dot X = F(X)$. If both branches give rise to approximately the same low-D trajectories, then we can say we have reduced our system. 
	
	Moreover, we need a nontrivial reduction that retains enough information of the statistics of the system in the high-dimensional space. For instance, the diagram always commutes for a single point projection, which of course erases all information of the original system. Thus we define a lift operator $\Psi$ that describes the statistic of the original system as a function of the reduced dynamics.  
	
	We will demonstrate the above described methods using the following system of individual modules:
	
	\begin{eqnarray}
	\dot{x}&=&w_{xx}x+w_{xy}y+f\left(x\right)\label{eq:module}\\
	\dot{y}&=&w_{yx}x+w_{yy}y+g\left(y\right)\nonumber 
	\end{eqnarray}

	\noindent The nonlinearities are captured by the analytic functions
	$f$ and $g$ with $f\left(0\right)=g\left(0\right)=Df\left(0\right)=Dg\left(0\right)=0$.
	
	Let $i\in\left\{ 1,2,\ldots,N\right\} $ be a parametrization of modules. We replicate the intramodule dynamics (\ref{eq:module}) for each of the $N$ modules and couple all modules through variable x with a linear connection matrix $M$ (and through layer y with a linear connection matrix $L$, but we will later consider $L=0$ ):
	
	\begin{eqnarray}
	\dot{x}_{i}&=&w_{xx}x_{i}+w_{xy}y_{i}+f\left(x_{i}\right)+\underset{j}{\sum}M_{ij}x_{j}\label{eq:2d_General}\\
	\dot{y}_{i}&=&w_{yx}x_{i}+w_{yy}y_{i}+g\left(x_{i}\right)+\underset{j}{\sum}L_{ij}y_{j}\nonumber 
	\end{eqnarray}

	A common way to study such systems is to define $M_{ij} = \lambda G_{ij}$, where $\lambda$ is a global coupling constant and $G_{ij}$ a graph connectivity matrix, whose elements are either 0 or 1, which can be used to describe the connections as occurring only along an underlying lattice. This is a standard setting in physics, where the interaction strength is controlled by a physical constant and therefore affects all extant interactions in parallel. This approach leads to a changing effective dimensionality of the dynamical system, because as $\lambda$ is varied, all dynamical eigenvalues of the system move in unison, and hence the number of modes which can potentially enter the dynamics increases with increasing $\lambda$. We will take a different approach.
	
		\begin{figure}[bth]
		\begin{centering}
			\includegraphics[scale=0.4]{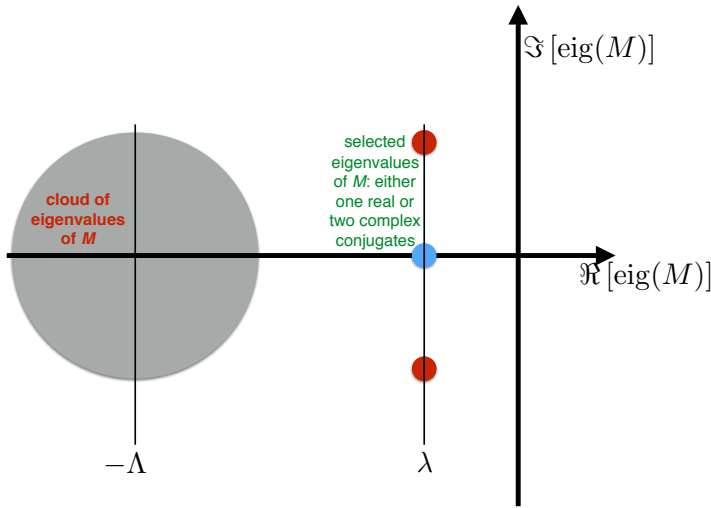}
			\par\end{centering}
		
		\begin{centering}
			\caption{\label{fig:eigenvalue_diagram} \textbf{The Eigenvalues of $M$ in the Complex Plane} The eigenvalues of $M$ in the complex plane. We shall construct $M$ so that most eigenvalues have strongly negative real parts, around $-d$, while either one real or two complex-conjugate eigenvalues with real part $\lambda$ will be allowed to approach the imaginary axis.}
			\par\end{centering}
		
	\end{figure}	
	
	We shall assume the matrix $M_{ij}$ is under control from slow homeostatic processes, and we shall allow either one single real eigenvalue or a couple of complex-conjugate eigenvalues to increase and approach the stability limits from the left (Fig. \ref{fig:eigenvalue_diagram}). The classical approach of doing this in the field of random neural networks is by first generating a base matrix $M_0$ with a given controlled spectrum. Following a long history of modeling studies, one can assume $M_0$ to be given by suitably-scaled i.i.d. Gaussian random variables \cite{Sompolinsky1998chaos}. Following this line of thought, in this work we shall either use $M_0$ as an $N \times N$ array of i.i.d Gaussians of variance $1/N$ (spectrum in the unit disk in the complex plane, matrix almost surely non-normal) or the antisymmetric component of said $M_0$, $M_a = (M_0 - M^t_0)/\sqrt{2}$ (purely imaginary spectrum, matrix is normal). The simplest way of then moving a single eigenvalue in the spectrum cloud independent of the rest of the eigenvalues and eigenvectors is the Brauer shift formula \cite{Brauer1952}:
	
	\begin{equation}
	M_{\mu}=M_{0}+(\lambda-\lambda_{0})({ \textbf{vw}})\label{eq:brauer}
	\end{equation}
	
	\noindent where $v$ and $w$ are the corresponding right and left eigenvectors.
	
	We can use this formula to simulate this system with rest of the eigenvalues damping at any value or for simplicity we can assume that the rest of the eigensystem has strongly negative eigenvalues. For large enough $d\in\mathbb{R}^{+}$, $M=M_{0/a}-d{ \textbf{I}}$ has the effect of moving the cloud of eigenvalues from a disk centered at zero to a disk centered at $-d$. We can then move a single real eigenvalue $\lambda\in\mathbb{R}$ (with $v$ and $w$ the corresponding right and left eigenvectors) to $\mu\in\mathbb{R}$ using the Brauer shift formula \cite{Brauer1952}: 
	
	\begin{equation}
	M_{\mu}=M_{0/a}-d{\textbf{I}}+(d+\mu-\lambda)({ \textbf{vw}})\label{eq:move_real}
	\end{equation}
	
	\noindent or a pair of complex conjugate eigenvalues $\lambda,\:\bar{\lambda}$
	to $\mu\pm im\left(\lambda\right)i$ :
	
	\begin{equation}
	M_{\mu}=M_{0/a}-d{ \textbf{I}}+(d+\mu-Re\left(\lambda\right))({ \textbf{vw}})+(d+\mu-Re\left(\bar{\lambda}\right))({ \mathbf{\bar{v}\bar{w}}})\label{eq:move_complex}
	\end{equation}.
	
	Antisymmetric matrices are normal. Thus consider $M$ to be given by eqn \ref{eq:move_real} with $M_a$, the antisymmetric component of a matrix of $N\times N$ i.i.d. Gausian variables of std $1/\sqrt{N}$; its eigenvalues are purely imaginary $\lambda = i \omega$ with imaginary component distributed like $p(\omega)\approx \sqrt{1-\omega^2/4}$. Consider one of its null eigenvalues (one is guaranteed to exist if $N$ is odd) and the corresponding  eigenvector $ v $. Then define as per Eq (\ref{eq:move_real})
	\begin{equation} 
	M_\lambda = M_a - dI + (d + \lambda) ( \mathbf{ v v^\intercal}) \label{eq:move_real_Lamda}
	\end{equation}
	where the left (dual) eigenvector $\mathbf{w}$ of $ \mathbf{v}$ is just it's transpose $ \mathbf{v}^\intercal$. All eigenvalues of $M$, except for the null one we chose, get moved to $-d$, while the null one now has value $\lambda$. 
	
	\noindent The Jacobian matrix at the fixed point $0$ is $J=\left[\begin{array}{c|c}
	M+w_{xx}I & w_{xy}I\\
	\\
	\hline \\
	w_{yx}I & L+w_{yy}I
	\end{array}\right]$.

	Both intermodule and intramodule connectivity and their interplay is essential in evaluating the eigenvalues - eigenvectors pairs of the Jacobian. A useful tool for this evaluation is the following lemma: 
	
	\noindent \underline{Lemma 1:} Let $W=\left[\begin{array}{c|c}
	M & w_{xy}I\\
	\\
	\hline \\
	w_{yx}I & L
	\end{array}\right]$ with $ML=LM$. Then :
	
	(a) $Wv=\lambda v\iff\lambda^{2}-\left(\lambda_{M}+\lambda_{L}\right)\lambda+\left(\lambda_{M}\lambda_{L}-w_{xy}w_{yx}\right)=0\,\;\&\,\;v=\begin{pmatrix}u\\
	\phi u
	\end{pmatrix},$ where $Mu=\lambda_{M}u,\,\;Nu=\lambda_{L}u$ and $\phi=\dfrac{\lambda-\lambda_{M}}{w_{xy}}=\dfrac{w_{yx}}{\lambda-\lambda_{L}}$
	. 
	
	\noindent \underline{Proof:} Trivial using that 2 matrices commute iff they are simultaneously diagonalizable.   $\blacksquare$

	The system as a whole, then, has a stability structure given by $2N$ eigenvalues that is a complex interplay of the structured column and the unstructured connections.

	
	\section{Example Network Used in This Paper}
	
	We now consider a particular case which we can solve in pretty good detail, so as to concentrate on what are the central figures and concepts in this approach. Since each individual unit is 2-D obviously it cannot do much more than have a number of fixed points and limit cycles. In what follows we will deal exclusively with the following $N$ coupled identical damped oscillators $\left(\ddot{x}_{i}+F\left(x_{i}\right)\dot{x}_{i}+x_{i}=\underset{j}{\sum}M_{ij}\dot{x}_{j}\right)$,
	i.e. for $i\in\left\{ 1,\ldots N\right\} $:
	\begin{eqnarray}
	\dot{x}_{i}&=&y_{i}+f\left(x_{i}\right)+\underset{j}{\sum}M_{ij}x_{j}\label{eq:leyerx}\\
	\dot{y}_{i}&=&-x_{i}\nonumber 
	\end{eqnarray}

	\noindent where $f\left(x\right)=$ $\underset{n\geq2}{\sum}c_{n}x^{n}$
	is analytic. 
	
	For example, in the special case $f(x) = -x^3$ we have $N$ coupled identical van der Pol oscillators:
	$$ \ddot x_i = -3\dot x_i x_i^2 - x_i + \sum_j M_{ij} \dot x_j $$
	
	We will show below we can evaluate the reduced dynamics for a Taylor series of $f$ order-by-order, so some of what follows concentrates on $f(x)=x^\alpha$. But the main issue is the structure of $M$, whether it's normal or non-normal; and whether we control a single real or two complex-conjugate eigenvalues. 
	
	As we shall see shortly, our procedure for obtaining the global dynamics from the micro-scopic dynamics involves both the right and the left eigenvectors of $M$. In general the matrix of left eigenvectors is the inverse of the matrix of right eigenvectors. However if $M$ is normal, then the inverse reduces to the complex conjugate making their relationship much simpler. Therefore we shall explore first normal matrices and then examine the additional issues raised by non-normality. Within normal matrices we shall, as anticipated in (Fig. \ref{fig:eigenvalue_diagram}), either control a single real eigenvalue or two complex-conjugate ones.
	
	We can evaluate the eigenvalue-eigenvector pairs of the Jacobian of the system from the corresponding eigenvalue-eigenvector pairs of $M$ using the following corollary:   
	
	\noindent \underline{Corollary 2:} Let $W=\left[\begin{array}{c|c}
	M & w_{xy}I\\
	\\
	\hline \\
	w_{yx}I & 0I
	\end{array}\right]$ with $w_{xy}w_{yx}<0$. Then:
	
	(a) $Wv=\lambda v\iff\lambda+\dfrac{-w_{xy}w_{yx}}{\lambda}=\lambda_{M}\,\;\&\,\;v=\begin{pmatrix}\lambda u\\
	w_{yx}u
	\end{pmatrix}$, where $Mu=\lambda_{M}u$
	
	(b)$sign\left(Re\left(\lambda\right)\right)=sign\left(Re\left(\lambda_{M}\right)\right)$
	
	(c)If $M$ has a zero eigenvalue , then $W$ has a pair of complex
	conjugate eigenvalues $\pm i\sqrt{-w_{xy}w_{yx}}$ 
	
	(d)If M has a pair of complex conjugate eigenvalues $\pm\alpha i$
	then $W$ has 2 pair of complex conjugate eigenvalues with 2 frequencies
	$\dfrac{\alpha\pm\sqrt{\alpha^{2}-4w_{xy}w_{yx}}}{2}$.
	
	\noindent \underline{Proof:} Trivia corollary of lemma 1 or see \cite{hoppensteadt2012weakly}
	Theorem 12.1 for a similar statement. $\blacksquare$

	 \section{Connectivity Matrix $M$}
	
	Let us assume that the elements of $M=\left(m_{i,j}\right)$ are independent randomly distributed gaussians of variance $\dfrac{1}{N}$ eg:  $m_{i,j} \  i.i.d. \  \sim\mathcal{N}\left(0,\dfrac{1}{N}\right)$ with spectrum in the unit disk in the complex plane\cite{Tao2008circular}.
	Let $\left\{ V_{,1},\ldots,V_{,N}\right\} $ be the basis in $\mathbb{C}^{N\times 1}\cong C^{\mathrm{N}}$
	of normalized right eigenvectors of $M\in\mathbb{R}^{N\times N}$ i.e. $MV=V\begin{pmatrix}\lambda_{1} &  & 0\\
	& \ddots\\
	0 &  & \lambda_{N}
	\end{pmatrix}$ , where $V=\begin{bmatrix}V_{,1}, & \cdots & ,V_{,N}\end{bmatrix}\in\mathbb{C}^{N\times N}$.
	Let $\left\{ W_{i,}=V_{,i}^{*}\right\} _{i\in\left\{ 1,\ldots N\right\} }$
	be the dual basis of left eigenvectors of $M$ in $\mathbb{R}^{1\times N}\cong\left(\mathbb{R}^{\mathrm{N}}\right)^{*}$
	(where $*:\mathbb{R^{\mathrm{N}}}\rightarrow\left(\mathbb{R}^{\mathrm{N}}\right)^{*}$
	is the dual map induced by the bilinear form $<V_{,i},V_{,j}>=\delta_{i,j}$)
	i.e. $W=\begin{bmatrix}W_{1,}\\
	\vdots\\
	W_{N,}
	\end{bmatrix}=V^{-1}$.
	
	For any $k\in\left\{ 1,\ldots,N\right\} $ and for any sequence of natural numbers of finite length $N$ : $p:\left\{ 1,\ldots,N\right\} \rightarrow\mathbb{N}$, and for a given ordered basis (in this case the right eigenvectors with a chosen ordering) let's define: 
	\begin{equation}
	\Gamma_{k}^{p}\coloneqq\underset{\varphi}{\sum}\left(V_{\varphi,k}\right)^{*}V_{\varphi,1}^{p\left(1\right)}\ldots V_{\varphi,N}^{p\left(N\right)}=\underset{\varphi}{\sum}W_{k,\varphi}V_{\varphi,1}^{p\left(1\right)}\ldots V_{\varphi,N}^{p\left(N\right)}=W_{k,}\times\underset{i}{\circledast}V_{,i}^{p\left(i\right)}\label{eq:gammas}
	\end{equation}
	where $\times$ is matrix multiplication and $\circledast$ is element-wise multiplication (and the powers $V_{,i}^{p\left(i\right)}$ are elementwise).
	
	\noindent In particular for the sequence  $\delta_{i,\bullet}$ defined by  $j \mapsto \delta_{i,j}$ (\ref{eq:gammas}) defines for $k=i$:
	
	\begin{equation}
	\Gamma_{i}^{k\delta_{i,\bullet}}\coloneqq\underset{\varphi}{\sum}\left(V_{\varphi,i}\right)^{*}V_{\varphi,i}^{k}=\underset{\varphi}{\sum}W_{i,\varphi}V_{\varphi,i}^{k}\label{eq:gammas2}
	\end{equation}
	
	As we will see these random variables play crucial role in describing the essential dynamics of the system. All the results in this part depend on $M$ only through the distribution of the $\Gamma_{k}^{p}$s and therefore can be extended to any connectivity type (for instance sparse matrices or only close neighbors connections) as long as we know the statistical properties of the $\Gamma_{k}^{p}$s. As we will see the limits: $$\underset{N\rightarrow\infty}{\lim}\left(N^{\frac{n-1}{2}}\dfrac{n!}{k_{1}!\cdots k_{N}!}\Gamma_{\mu}^{\left(k_{1},\ldots,k_{N}\right)}\right)$$ where $n=k_{1}+\ldots+k_{N}$ are crucial, and in particular the limits: $\underset{N\rightarrow\infty}{\lim}\left(N^{\frac{k-1}{2}}\Gamma_{i}^{k\delta_{i,\bullet}}\right)$. 
	
	\section{Calculate $\Gamma$'s for Normal Matrix with Antisymmetric Gaussian i.i.d. Entries}

	In the case of normal matrices (\ref{eq:gammas2}) can be written us:
	
	\begin{equation}
	\Gamma_{i}^{k\delta_{i,\bullet}}=\underset{\varphi}{\sum}V_{\varphi,i}^{k+1}\label{eq:gammas2normal}
	\end{equation}
	
	In the case of antisymmetric matrices we numerically observe that the elements of the eigenvectors approximate Gaussians of variance $\dfrac{1}{N}$ as $N\rightarrow \infty$ with:
	
	\begin{equation} 
	\underset{N\rightarrow\infty}{\lim}\left(N^{\frac{k-1}{2}}\Gamma_{i}^{k\delta_{i,\bullet}}\right)  = 
	\begin{cases}
	\ k !! \   & k \text{ odd, } \\
	\ 0 & k \text{ even.}
	\end{cases}
	\end{equation}
	
	\noindent given by the moments of Gaussian, where $k!!$ is the double factorial $k!!=1\cdot3\cdots(k-2)\cdot k$. Similarly :
	
	\begin{equation}
	\underset{N\rightarrow\infty}{\lim}\left(N^{\frac{n-1}{2}}\dfrac{n!}{k_{1}!\cdots k_{N}!}\Gamma_{i}^{\left(k_{1},\ldots,k_{N}\right)}\right)=\dfrac{n!}{k_{1}!\cdots k_{N}!}\left(k_{1}-1\right)!!\cdots k_{i}!!\cdots\left(k_{N}-1\right)!! = n!!\begin{pmatrix}\frac{n-1}{2}\\
	\frac{k_{1}}{2}\cdots\frac{k_{i}-1}{2}\cdots\frac{k_{N}}{2}
	\end{pmatrix} 
	\end{equation}
	
	\noindent if $k_{i}$ and $n$ are odd and $k_{j}$ for $j\neq i$ are even. Else $=0$. So the  $\Gamma$'s are given by co-moments of independent Gaussians.

	Let us notice that it is known that for random symmetric matrices (Gaussian orthogonal ensemble) eigenvectors are uniformly distributed in the sphere and independent of the eigenvalues (Corollary 2.5.4 in \cite{Anderson2010}). This is implied by the invariance of the law of X under arbitrary orthogonal transformations as in our case. From \cite{O'Rourke2016Eigenvectors} we know that if $v$ is uniformly distributed on the sphere then $\lVert{v_{i}}\rVert^{2}\sim Beta\left(\dfrac{1}{2},\dfrac{N-1}{2}\right)$ and $\underset{N\rightarrow\infty}{\lim}\lVert{N^{\frac{p}{2}-1}\lVert{v}\rVert_{l_{p}}^{p}-E\lvert{N\left(0,1\right)}\rvert}^{p}\rVert=0$ almost surely.

	\section{Linear Projection onto the Stable and the Center Spaces - Definitions}
	Let $X_{\lambda_{k}}=\dfrac{1}{\sqrt{N}}W_{k,}\times\begin{pmatrix}x_{1}\\
	\vdots\\
	x_{n}
	\end{pmatrix}=\dfrac{1}{\sqrt{N}}\underset{i}{\sum}W_{k,i}x_{i}$ be the normalized (so that $X_{\lambda_{k}}$ is of order 1, because ${ x \cdot x} \approx N$ while if $M$ is normal $\ W_{k,} \cdot W_{k,} = 1$ so that $ { W_{k,} \cdot x } \approx \sqrt{N}$ generically.)   coordinates of activity in the basis of right eigenvalues, 
	so :

	\begin{equation}
	x=\begin{pmatrix}x_{1}\\
	\vdots\\
	x_{N}
	\end{pmatrix}=\sqrt{N}V\begin{pmatrix}X_{\lambda_{1}}\\
	\vdots\\
	X_{\lambda_{N}}
	\end{pmatrix}=\sqrt{N}\underset{i}{\sum}X_{\lambda_i}V_{,i}\label{eq:eigbaseTobase}
	\end{equation}
	
	\noindent where $V_{,i}$, for $i=1,...,N$ are the right eigenvectors, columns of the matrix $V$.
	
	\noindent Applying $\dfrac{1}{\sqrt{N}}W$ to equations (\ref{eq:leyerx}) we have, for an eigenvalue $\mu\in spec\left(M\right)$
	:
	
	\begin{eqnarray}
	\dot{X}_{\mu}&=&\dfrac{1}{\sqrt{N}}\underset{i}{\sum}W_{\mu,i}f\left(x_{i}\right)+Y_{\mu}+\mu X\label{eq:CoarseGrain1}\\
	\dot{Y}_{\mu}&=&-X_{\mu}\nonumber 
	\end{eqnarray}

	Equation (\ref{eq:CoarseGrain1}) will be essential in defining the operator $\Phi$ as seen in (Fig. \ref{fig:schematics2}). Equation (\ref{eq:eigbaseTobase}) will be essential in defining the operator $\Psi$.
	
	\noindent The Hartman\textendash Grobman theorem only holds for hyperbolic equilibriums, so as eigenvalues move close to the critical boundary, the nonlinearities:
	
	\begin{equation}
	F_{\mu}\coloneqq\dfrac{1}{\sqrt{N}}\underset{i}{\sum}W_{\mu,i}f\left(x_{i}\right)=\dfrac{1}{\sqrt{N}}\underset{i}{\sum}W_{\mu,i}f\left(\sqrt{N}\underset{k}{\sum}V_{i,k}X_{\lambda_{k}}\right)\label{eq:NonLinearities_CoarseGrain}
	\end{equation}

	\noindent became essential in describing the dynamics.
	
	Equations (\ref{eq:leyerx}) (similarly for the general equations (\ref{eq:2d_General})) can be written in the basis of right eigenvectors corresponding to the eigenvalues, for $\mu\in spec\left(M\right)$:
	
	\begin{flushright}
		\begin{eqnarray}
		& \dot{X}_{\mu}=F_{\mu}\left(X_{\lambda_{1}},\ldots,X_{\lambda_{N}}\right)+Y_{\mu}+\mu X_{\mu}\label{eq:System_Base_Change}\\
		& \dot{Y}_{\mu}=-X_{\mu}\nonumber 
		\end{eqnarray}
		
		\par\end{flushright}
	
	\noindent where 
	\begin{equation}
	F_{\mu_{i}}\left(X_{\lambda_{1}},\ldots,X_{\lambda_{N}}\right)=\underset{n\geq2}{\sum}c_{n}\sum_{k_{1}+\ldots+k_{N}=n}\left(\dfrac{n!}{k_{1}!\cdots k_{N}!}{\sqrt{N}}^{\left(n-1\right)}\Gamma_{\mu}^{\left(k_{1},\ldots,k_{N}\right)}\right)X_{\lambda_{1}}^{k_{1}}\cdots X_{\lambda_{N}}^{k_{N}}\label{eq:NonLinearities_Base_Change}
	\end{equation}
	
	\section{Naive Approach - Single Real Eigenvalue, Normal Matrix}
	
	Let us consider the system (\ref{eq:leyerx}) as a single real eigenvalue of M, $\lambda=\lambda_{1}\in\mathbb{R}$, crosses the imaginary axis (\ref{eq:move_real}) (the critical boundary for this case, look at corollary 2) becoming slightly positive, while the rest of the eigenvalues $\mu\in\left\{ \lambda_{2},\ldots,\lambda_{N}\right\} $ have negative real part. The naive approximation would be to consider
	$X_{\mu}\approx0$ for $\mu\in spec\left(M\right)\setminus\left\{ \lambda\right\} $ and we can reduce the $2N$-dimensional system (\ref{eq:leyerx})
	into the 2-dimensional system :

	\begin{eqnarray}
	\dot X  &=& F(x) + Y + \mu X  \label{eqmacro}\\
	\dot Y  &=& -X\nonumber
	\end{eqnarray}

	Where $X\coloneqq X_{\lambda}$ , $Y \coloneqq Y_{\lambda}$ and $F(x)$ is given by (\ref{eq:NonLinearities_CoarseGrain}):
	
	\begin{equation}
	F(x)\coloneqq F_{\lambda}(x)=\dfrac{1}{\sqrt{N}}\underset{i}{\sum}W_{\lambda,i}f\left(x_{i}\right)=\dfrac{1}{\sqrt{N}}\underset{i}{\sum}W_{\lambda,i}f\left(\sqrt{N}V_{i,1}X+\sqrt{N}\underset{\mu\in spec\left(M\right)\setminus\left\{ \lambda\right\}}{\sum}V_{i,\mu}X_{\mu}\right)\label{eq:NonLinearities_CoarseGrain1}
	\end{equation}
	
	Evaluation of the $\underset{i}{\sum}W_{\lambda,i}f\left(x_{i}\right)$ term is the key; assuming conversely that the projection of $x_i$ onto the complement of $V_{,1}$ (space perpendicular to $V_{,1}$ in the normal case) is small (an assumption we shall revisit in detail below) we'd get an ansatz which is our first, naive approximation of the lift operator $\Psi$, which gives the coordinates in the full space as a function of the coordinates on the center manifold ${ x} = \Psi[X]$:
	\begin{eqnarray} 
	x_i &\approx&   \sqrt{N} X V_{i,1} \label{eq:Psi_naive} \\
	y_i &\approx&  \sqrt{N} Y V_{i,1} \nonumber
	\end{eqnarray}
	
	Then (\ref{eqmacro}) can be written us : 
	
	\begin{eqnarray}
	\dot X  &=& F(X) + Y + \mu X  \label{eqmacro1}\\
	\dot Y  &=& -X\nonumber
	\end{eqnarray}
	
	Where :
	
	\begin{equation}
	F(X)=\dfrac{1}{\sqrt{N}}\underset{i}{\sum}W_{\lambda,i}f\left(\sqrt{N}V_{i,1}X\right)\label{eq:NonLinearities_CoarseGrain2}
	\end{equation}

	\noindent The coarse-graining operator $\Phi$, which transforms the original full dimensional dynamics to the coarsened, low-dimensional dynamics is defined by $$F \equiv \Phi[ f ]$$

	In other words the $2N$ dimensional equation (\ref{eq:CoarseGrain1}) is reduced to the 2 dimensional equation: 
	
	\begin{eqnarray}
	\dot{X}_{\lambda}&=&\underset{k\geq2}{\sum}\left(c_{k}{\sqrt{N}}^{\left(k-1\right)}\Gamma_{\lambda}^{k\delta_{1,\bullet}}\right)X_{\lambda}^{k}+Y_{\lambda}+\lambda X_{\lambda}\label{eq:layerx_CoarseGrain_1}\\
	\dot{Y}_{\lambda}&=&-X_{\lambda}\nonumber 
	\end{eqnarray}

	We see that the dynamics of $X_{\lambda}$ is determined by the random variables: $\sqrt{N}^{(k-1)}\Gamma_{\lambda}^{k\delta_{1,\bullet}}=\sqrt{N}^{(k-1)}\underset{\varphi=1}{\overset{N}{\sum}}W_{1,\varphi}V_{\varphi,1}^{k}$ where $V_{\cdotp,1}$ is normalized so that $\underset{\varphi}{\sum}V_{\varphi,1}^{2}=1$ :

	\noindent In the case of normal matrices and $f\left(x\right)=$ $\underset{n\geq2}{\sum}c_{n}x^{n}$ analytic we have: 
	\begin{equation}
	F\left(X\right)=\dfrac{1}{\sqrt{N}}\underset{\varphi=1}{\overset{N}{\sum}}W_{1,\varphi}f\left(V_{\varphi,1}\sqrt{N}X\right)=\underset{k\geq1}{\sum}\left(c_{2k+1}\left(2k+1\right)!!\right)X^{2k+1}
	\end{equation}. 
	
	The first thing to notice is that $\Phi$ is linear in $f$, meaning we can try to apply it order-by-order in a Taylor expansion of $f$. Applying Eq~(\ref{eq:NonLinearities_CoarseGrain2}) to individual integer powers yields, in the large $N$ limit, an evaluation of the moments of a Gaussian distribution: 
	\begin{equation} 
	\Phi[ x^\alpha ] = X^\alpha N^{{\alpha - 1}\over 2} \sum_i e_i^{\alpha+1}  = 
	\begin{cases}
	\ \alpha !! \  X^\alpha & \alpha \text{ odd, } \\
	\qquad 0 & \alpha \text{ even.}
	\end{cases}\label{eq:Naive_phi}
	\end{equation}
	The most important property to note is that the result is  independent of which eigenvector $ v$ we chose, because under these circumstances the operator is self-averaging. Then, the power-law is unchanged except for a prefactor: the cubic $X^3$ nonlinearity gets renormalized by a factor of 3, $X^5$ by 15, $X^7$ by 105. Curiously this operator {\em destroys all even nonlinearities}. As we shall see later the even nonlinearities get absorbed in the overall shape of the center manifold, which we have assumed in the ansatz Eq \eqref{eq:Psi_naive} to be linear in $x_i$. However the naive reduced equations Eq (\ref{eqmacro}) are not affected. 
	
\begin{figure}[]
		\begin{centering}
			\includegraphics[scale=0.6]{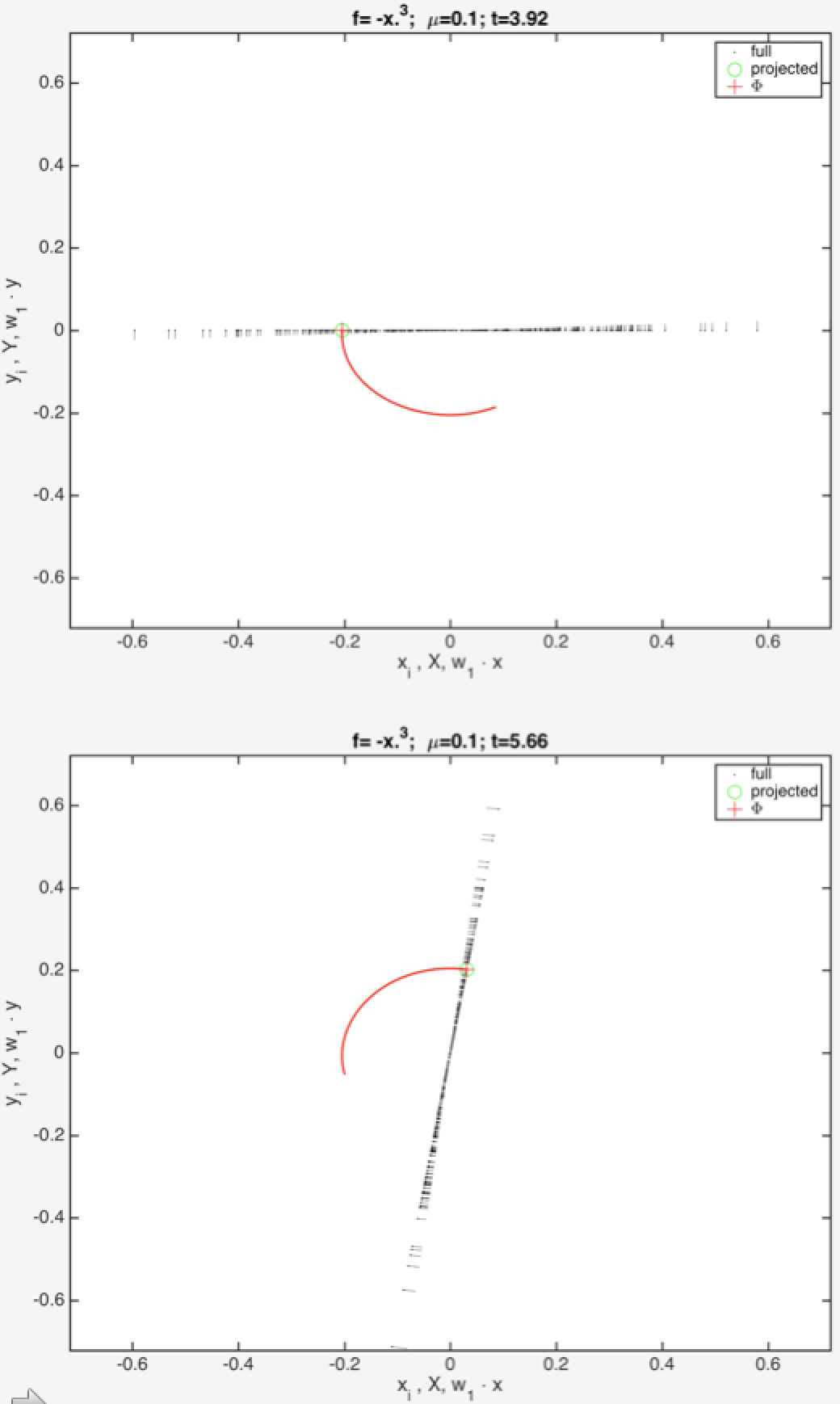}
			\par\end{centering}
		
		\begin{centering}
			\caption{\label{fig:naive1} \textbf{Naive $\Phi$ - Linear $\Psi$} Randomly coupled Van der Pol oscillators given by Eq (\ref{eq:leyerx}) with $f(x)=-x^3$ and for random matrix $M$ defined using Eq \eqref{eq:move_real_Lamda} with $N=351$ and $d=-30$. We let $\lambda=0.1$ cross the real line by a small value.}
			\par\end{centering}
		
	\end{figure}

	We first simulate randomly coupled Van der Pol oscillators given by Eq (\ref{eq:leyerx}) with $f(x)=-x^3$ and for random matrix $M$ defined using Eq \eqref{eq:move_real_Lamda} with $N=351$ and $d=-30$. We let $\lambda=0.1$ cross the real line by a small value. We numerically integrate the whole system for random initial values and let it settle after a transient on the limit cycle. In movie: \url{https://drive.google.com/open?id=14ku9fZlCUAW_mzEHIvezkRw4G6fbAAix} we plot at each time frame $t$ the values $\{(x_i(t),y_i(t))\arrowvert i\in\{1,\ldots,N\}$ as black dots. Moreover for each $i$ we plot as a thin black line the orbit of $(x_i,y_i)$ for a short time past interval $(t-\Delta t,t)$ to visualize the orbits. We plot with a green cycle $(X_\lambda,Y_\lambda)$ the linear projection of $x$ and $y$ in the critical eigenvector and with a red cross the simulation of Eq \eqref{eqmacro1} with $F(X)=-3X^3$ (and initial condition the linear projections on the critical mode of the initial conditions used for the high dimensional Eq (\ref{eq:leyerx})). Snapshots of the simulation movie are shown in (Fig \ref{fig:naive1}). As in the case of individual orbits we also plot part of the past orbits for visualization with green and red colors respectable. As we can see in this case diagram (Fig. \ref{fig:schematics2}) using our naive definition of our $\Phi$ operator Eq \eqref{eq:Naive_phi}. Moreover we see that the cubic nonlinearity doesn't survive in the lift  $\Psi$ a statement that we shall see theoretically shortly.

\begin{figure}[]
		\begin{centering}
			\includegraphics[scale=0.6]{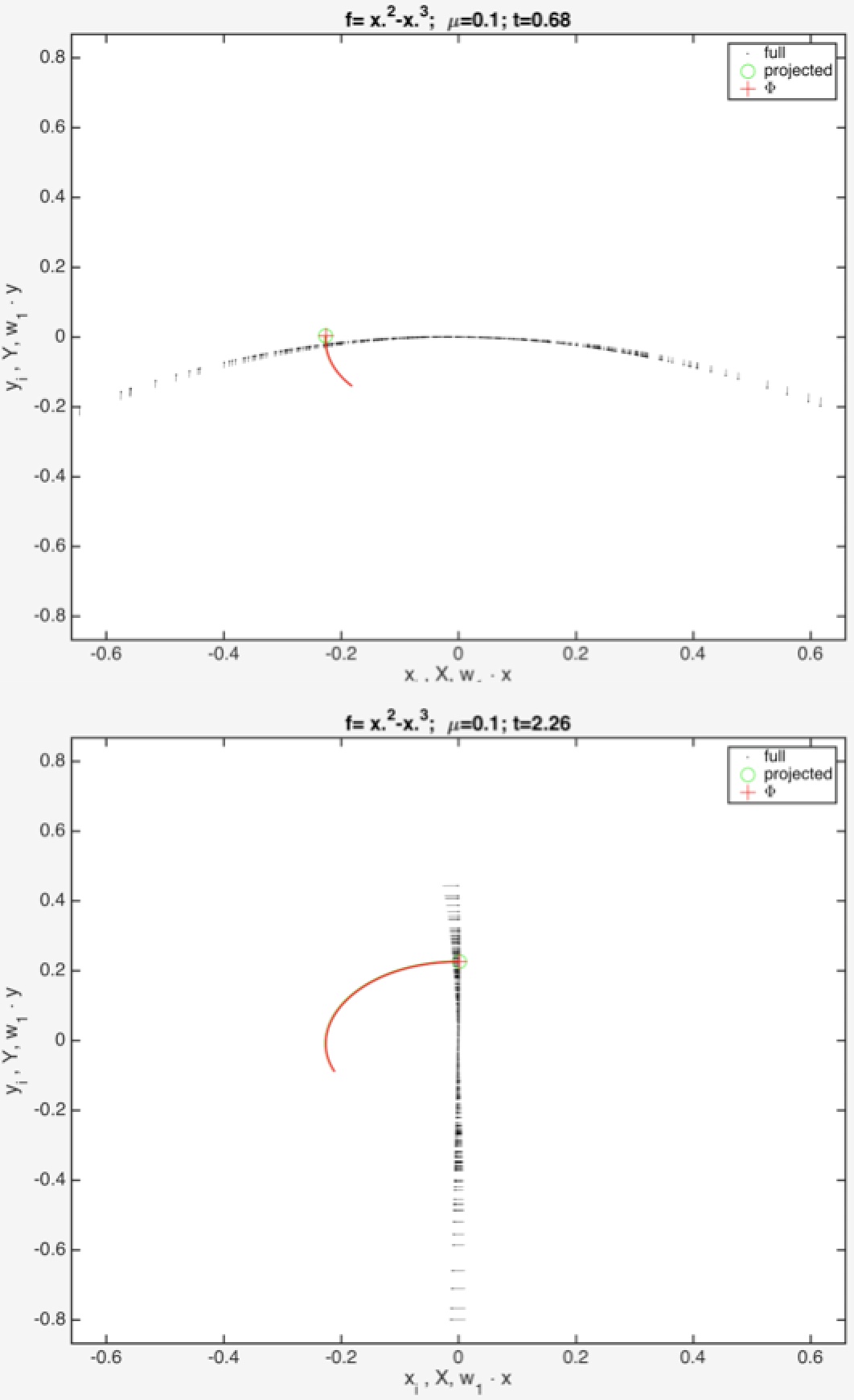}
			\par\end{centering}
		
		\begin{centering}
			\caption{\label{fig:naive2} \textbf{Naive $\Phi$ - Nonlinear $\Psi$} Randomly coupled Van der Pol oscillators given by Eq (\ref{eq:leyerx}) with $f(x)=x^2-x^3$ and for random matrix $M$ defined using Eq \eqref{eq:move_real_Lamda} with $N=351$ and $\Lambda=-30$. We let $\lambda=0.1$ cross the real line by a small value.}
			\par\end{centering}
		
	\end{figure}

	In movie: \url{https://drive.google.com/open?id=1Ldn_YvlXaCSRc9L3dS_5smrh_39OqzOP} we can see the corresponding simulation for $f(x)=x^2-x^3$. Snapshots of the simulation movie are shown in (Fig. \ref{fig:naive2}). As we can see the added even nonlinearity doesn't survive in the $\Phi$ operator and our approximation is still accurate (as we will see later a more accurate statement would be that the even nonlinearity doesn't survive as an even nonlinearity but can have small "leak" into higher order odd nonlinearities). But in this case it becomes apparent that the naive linear approximation of $\Psi$ Eq \eqref{eq:Psi_naive} fails and the even nonlinearity survives in the operator. We will see that theoretically shortly.  
	
	\section{Failure of Above}
	
	The above approximation is studied in \cite{Moirogiannis2017} and assumes that the center space $E^{c}=\mathbb{R}\begin{pmatrix}v_{1}\\
	0
	\end{pmatrix}\oplus\mathbb{R}\begin{pmatrix}0\\
	v_{1}
	\end{pmatrix}$ is invariant which is only true up to linear approximation. In fact the nonlinearities of the invariant center manifolds tangent to $E_{c}$ can change the qualitative behavior of the reduced to the manifold dynamics and there are examples that show that the above naive linear projection fails to capture the reduced dynamics \cite{Guckenheimer2013}. In our model particularly equations (\ref{eq:layerx_CoarseGrain_1}) predict that the Hopf bifurcation
	is supercritical iff $c_{3}<0$ independently of the value of $c_{2}$ but as numeric simulations show and as we will prove later the bifurcation is supercritical for large enough $\left|c_{2}\right|$. 
	
\begin{figure}[]
		\begin{centering}
			\includegraphics[scale=0.7]{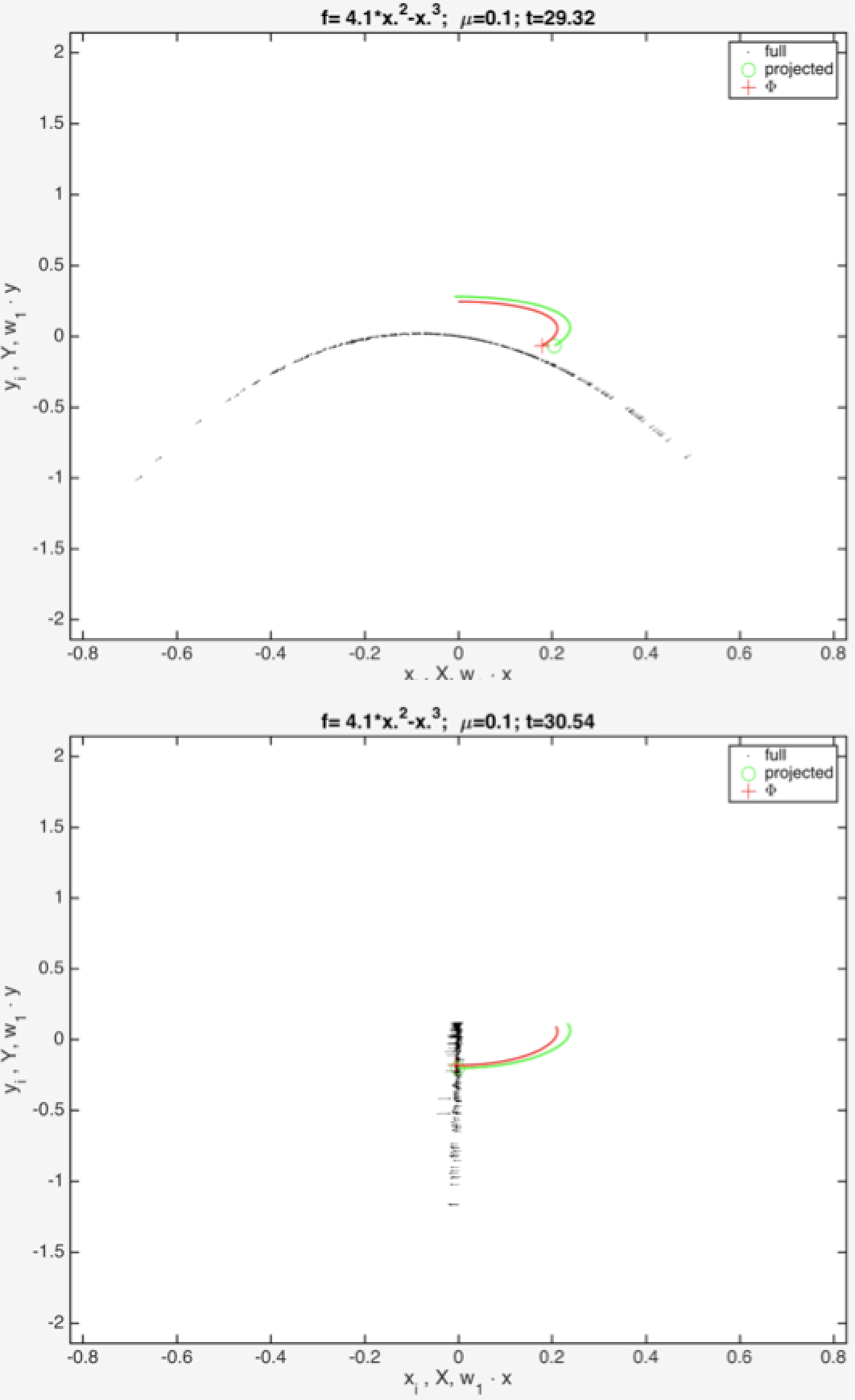}
			\par\end{centering}
		
		\begin{centering}
			\caption{\label{fig:naive_phi_fail} \textbf{Naive $\Phi$ Fail} Randomly coupled Van der Pol oscillators given by Eq (\ref{eq:leyerx}) with $f(x)=4.1*x^2-x^3$ and for random matrix $M$ defined using Eq \eqref{eq:move_real_Lamda} with $N=351$ and $\Lambda=-30$. We let $\lambda=0.1$ cross the real line by a small value.}
			\par\end{centering}
		
	\end{figure}

	We can see the failure of our naive $\Phi$ approximation Eq \eqref{eq:Naive_phi} by simulating a system as in the last movie with a higher order of even nonlinearity $f(x)=4.1*x^2-x^3$. We can see this in movie: \url{https://drive.google.com/open?id=1FpkxPm5kRCjHi_TSl_5Mrk3Jn6HQgAlr}. Snapshots of the simulation movie are shown in (Fig. \ref{fig:naive_phi_fail}).

	\section{Center Manifold Reduction}
	
	We will follow similar notation and approach as in \cite{hoppensteadt2012weakly} ,
	in particular look at section 12.7 for the 3rd order approximation
	of $\psi$ for $f\left(x\right)=x^{2}$ and 2nd order approximation
	of $\psi$ for $f\left(x\right)=\pm x^{3}$ in the case of symmetric
	interaction matrix. Here we care about more general forms of matrices and for any abstract nonlinearity. For more details on center manifolds look at \cite{Ioos1999} or for formal proofs at \cite{Carr2012}.
	
	We can rewrite equation \eqref{eq:leyerx} as: 
	
	\noindent $\eqref{eq:leyerx}\Leftrightarrow\quad$ : $\dot{Z}=G\left(Z\right)=WZ+F\left(Z\right)$

	\noindent where $Z=\left(\begin{array}{c}
	x_{1}\\
	\vdots\\
	x_{N}\\
	y_{1}\\
	\vdots\\
	y_{N}
	\end{array}\right)\in\mathbb{R}^{2N},\quad W=\left(\begin{array}{c|c}
	M & I\\
	\\
	\hline \\
	-I & 0
	\end{array}\right),\quad F\left(Z\right)=\left(\begin{array}{c}
	f\left(x_{1}\right)\\
	\vdots\\
	f\left(x_{N}\right)\\
	0\\
	\vdots\\
	0
	\end{array}\right)$ . 
	
	Let us assume that a single real eigenvalue of M, $\lambda=\lambda_{1}\in\mathbb{R}$ (with corresponding right eigenvector $v_{1}$) crosses the imaginary axis. Corollary 2 implies that the corresponding eigenvalues of the whole system are $\pm i$ with corresponding eigenvectors $\begin{pmatrix} \pm iv_{1}\\
	-v_{1}
	\end{pmatrix}$. Let $E^{c}=\mathbb{R}\begin{pmatrix}v_{1}\\
	0
	\end{pmatrix}\oplus\mathbb{R}\begin{pmatrix}0\\
	v_{1}
	\end{pmatrix}=\left\{ x\begin{pmatrix}v_{1}\\
	0
	\end{pmatrix}+y\begin{pmatrix}0\\
	v_{1}
	\end{pmatrix}\mid x,y\in\mathbb{R}\right\} $ be the center space, $E^{s}=E^{\lambda_{2}}\oplus\ldots\oplus E^{\lambda_{N}}=$$\mathbb{R}\left(\begin{array}{c}
	\upsilon_{2}\\
	0
	\end{array}\right)\oplus\mathbb{R}\begin{pmatrix}0\\
	v_{2}
	\end{pmatrix}\oplus\ldots\oplus\mathbb{R}\left(\begin{array}{c}
	\upsilon_{N}\\
	0
	\end{array}\right)\oplus\mathbb{R}\begin{pmatrix}0\\
	v_{N}
	\end{pmatrix}$ be the stable space, $\Pi_{c}=\begin{pmatrix}v_{1}w_{1} & 0\\
	0 & v_{1}w_{1}
	\end{pmatrix}:\mathbb{R}^{2N}=E^{c}\oplus E^{s}\rightarrow E^{c}$ , $\Pi_{\lambda_{k}}=\begin{pmatrix}v_{k}w_{k} & 0\\
	0 & v_{k}w_{k}
	\end{pmatrix}:\mathbb{R}^{2N}\rightarrow E^{\lambda_{k}}$ and $\Pi_{s}=\Pi_{\lambda_{2}}+\ldots+\Pi_{\lambda_{k}}:\mathbb{R}^{2N}=E^{c}\oplus E^{s}\rightarrow E^{s}$ the defined by the basis projections. The above projections commute with $W$ since they have the same eigenspaces. 
	
	The center manifold theorem ensures that there is a local mapping $\psi:E^{c}\rightarrow E^{s}$ with 
	\[
	\psi(0)=0\quad and\quad D\psi(0)=0
	\]

	\noindent such that the manifold $\mathcal{M}$ defined by:
	\[
	\mathcal{M}=\left\{ u+\psi\left(u\right)\mid u\in E_{c}\right\} =\left\{ x\begin{pmatrix}v_{1}\\
	0
	\end{pmatrix}+y\begin{pmatrix}0\\
	v_{1}
	\end{pmatrix}+\psi\left(x,y\right)\mid x,y\in\mathbb{R}\right\} 
	\]

	\noindent is invariant and locally attractive (emergence theorem).
	The dynamics reduced on $\mathcal{M}$ are locally:
	
	\begin{equation}
	\dot{u}=\Pi_{c}G\left(u+\psi\left(u\right)\right)=Wu+\Pi_{c}F\left(u+\psi\left(u\right)\right)\label{eq:reduction}
	\end{equation}

\begin{figure}[]
		\begin{centering}
			\includegraphics[scale=0.5]{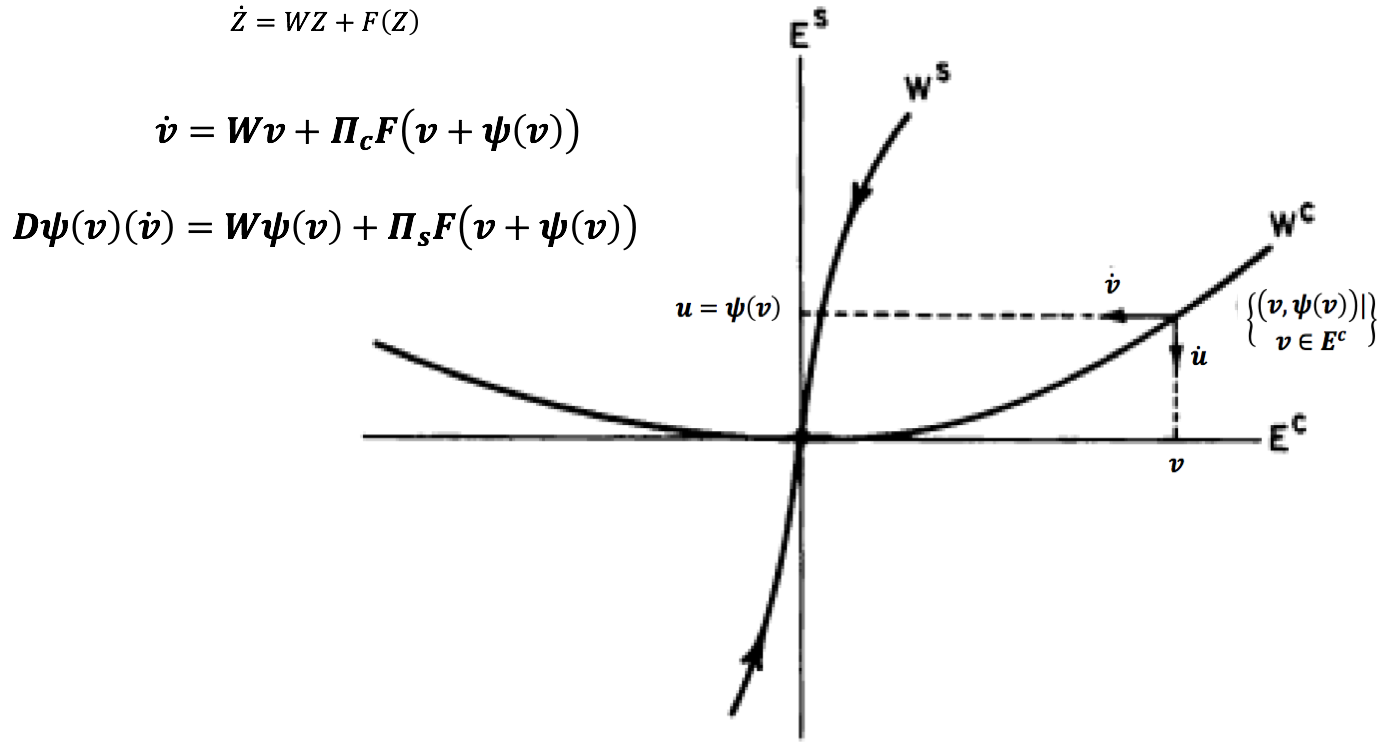}
			\par\end{centering}
		
		\begin{centering}
			\caption{\label{fig:center_manifold} \textbf{Center Manifold Reduction} A schematic of center manifold and the invariance equations.}
			\par\end{centering}
		
	\end{figure}

	\noindent A schematic of a center manifold and the invariance equations can be seen in (Fig. \ref{fig:center_manifold}).  Using the approximation theorem, we can approximate $\psi$ to any order from the equation:
	
	\begin{equation}
	D\psi(u)\left(Wu+\Pi_{c}F\left(u+\psi\left(u\right)\right)\right)=D\psi(u)\dot{u}=\Pi_{s}G\left(u+\psi\left(u\right)\right)=W\psi\left(u\right)+\Pi_{s}F\left(u+\psi\left(u\right)\right)\label{eq:approx}
	\end{equation}

	\noindent More precisely, if $Dh(u)\left(Wu+\Pi_{c}F\left(u+h\left(u\right)\right)\right)-Wh\left(u\right)+\Pi_{s}F\left(u+h\left(u\right)\right)=O\left(\left|u\right|^{k}\right)$
	then $\psi\left(u\right)=h\left(u\right)+O\left(\left|u\right|^{k}\right)$.
	
	\noindent $\psi$ approximated by (\ref{eq:approx})  is used to define the lift operator $\Psi_{f}$ and equation
	(\ref{eq:reduction}) defines the reduction operator $\Phi$.

	This gives us an algorithm of calculating orders of increasingly degree in the Taylor expansions of both reduced dynamics and the center manifold embedding by iteratively interchanging  equations (\ref{eq:reduction}) : $\dot{u}=Wu+\Pi_{c}F\left(u+\psi\left(u\right)\right)$ and (\ref{eq:approx}) : $D\psi(u)\dot{u}=W\psi\left(u\right)+\Pi_{s}F\left(u+\psi\left(u\right)\right)$ correspondingly. However this algorithm can be complicated. As the dimension of the center manifold increases (or order of desired approximation increases) the calculations become intractable (the number of unknown parameters and resulting equations grows like $\begin{pmatrix}dimension+order-1\\ order\end{pmatrix},$ the number of parameters to keep track of $\sim$$\underset{r<order}{\sum}\begin{pmatrix}dimension+r-1\\r\end{pmatrix}$ + all combinations with repetitions of those up to the specified $order$). So using symbolic algorithms is essential \cite{Freire1988Algorithm} .
	
	\section{$n$ th Order Approximation of $\Psi$ for $f\left(x\right)=c_{n}x^{n}+O\left(x^{n+1}\right)$ }
	
	Even though those questions are in general complicated, when estimating the first nonlinear order of approximation for the center manifold embedding and reduced dynamics on it of \eqref{eq:leyerx} with $f\left(x\right)=c_{n}x^{n}+O\left(x^{n+1}\right)$ the corresponding equations are very simple.
	
	Let $\psi\left(\begin{pmatrix}v_{1}\\
	0
	\end{pmatrix}x+\begin{pmatrix}0\\
	v_{1}
	\end{pmatrix}y\right)=\underset{k=2}{\overset{N}{\sum}}\left(\begin{array}{c}
	\left(\underset{i=0}{\overset{n}{\sum}}\alpha_{i}^{k}x^{n-i}y^{i}\right)V_{,k}\\
	\left(\underset{i=0}{\overset{n}{\sum}}\beta_{i}^{k}x^{n-i}y^{i}\right)V_{,k}
	\end{array}\right)+O\left(\left(\sqrt{x^{2}+y^{2}}\right)^{n}\right)$ , and for convenience let $\alpha_{-1}=\alpha_{n+1}=\beta_{-1}=\beta_{n+1}=0$.
	
	\noindent Then $D\psi\left(u\right)\dot{u}=\underset{k=2}{\overset{N}{\sum}}\left(\begin{array}{c}
	\left(\underset{i=0}{\overset{n}{\sum}}\left(\left(n-i+1\right)\alpha_{i-1}^{k}-\left(i+1\right)\alpha_{i+1}^{k}\right)x^{n-i}y^{i}\right)V_{,k}\\
	\left(\underset{i=0}{\overset{n}{\sum}}\left(\left(n-i+1\right)\beta_{i-1}^{k}-\left(i+1\right)\beta_{i+1}^{k}\right)x^{n-i}y^{i}\right)V_{,k}
	\end{array}\right)+O\left(\left(\sqrt{x^{2}+y^{2}}\right)^{n}\right)$ ,
	
	$W\psi\left(u\right)=\underset{k=2}{\overset{N}{\sum}}\left(\begin{array}{c}
	\left(\underset{i=0}{\overset{n}{\sum}}\left(\lambda_{k}\alpha_{i}^{k}+\beta_{i}^{k}\right)x^{n-i}y^{i}\right)V_{,k}\\
	\left(\underset{i=0}{\overset{n}{\sum}}-\alpha_{i}^{k}x^{n-i}y^{i}\right)V_{,k}
	\end{array}\right)+O\left(\left(\sqrt{x^{2}+y^{2}}\right)^{n}\right)$ ,
	
	$\Pi_{s}F\left(u+\Psi\left(u\right)\right)=c_{n}\Pi_{s}\left(\begin{array}{c}
	\left(V_{1,1}x+O\left(n+1\right)\right)^{n}\\
	\vdots\\
	\left(V_{N,1}x+O\left(n+1\right)\right)^{n}\\
	0\\
	\vdots\\
	0
	\end{array}\right)=c_{n}\Pi_{s}\left(\begin{array}{c}
	\left(V_{1,1}x\right)^{n}\\
	\vdots\\
	\left(V_{N,1}x\right)^{n}\\
	0\\
	\vdots\\
	0
	\end{array}\right)+O\left(n+1\right)=c_{n}x^{n}\underset{k=2}{\overset{N}{\sum}}V_{,k}\varGamma_{k}^{n\delta_{1,\bullet}}+O\left(n+1\right)$. 
	
	\noindent Projecting equation (\ref{eq:approx}) onto the stable modes using
	$\Pi_{\lambda_{k}}$ for $k\geq2$ gives us the equations: 
	
	\begin{eqnarray}
	\begin{pmatrix}\alpha_{0}^{k}\\
	\vdots\\
	\alpha_{n}^{k}
	\end{pmatrix}&=&A\begin{pmatrix}\beta_{0}^{k}\\
	\vdots\\
	\beta_{n}^{k}
	\end{pmatrix}\label{eq:Center_Manifold_Parameters}\\
	\left(-\lambda_{k}id_{n}-A\right)\begin{pmatrix}\alpha_{0}^{k}\\
	\vdots\\
	\alpha_{n}^{k}
	\end{pmatrix}-\begin{pmatrix}\beta_{0}^{k}\\
	\vdots\\
	\beta_{n}^{k}
	\end{pmatrix}&=&\begin{pmatrix}c_{n}\varGamma_{k}^{n\delta_{1,\bullet}}\\
	\vdots\\
	0
	\end{pmatrix}\nonumber 
	\end{eqnarray}

	\noindent where $A=\begin{pmatrix}0 & 1\\
	-n & 0 & 2\\
	& -(n-1) & 0 & \ddots\\
	&  & \ddots & \ddots & \ddots\\
	&  &  & \ddots & 0 & n-1\\
	&  &  &  & -2 & 0 & n\\
	&  &  &  &  & -1 & 0
	\end{pmatrix}$ . 
	
	If $f\left(x\right)=c_{2}x^{2}+O\left(x^{3}\right)$ then after scaling normalization we get:
	\begin{equation}
	\Psi\left(\begin{pmatrix}v_{1}\\
	0
	\end{pmatrix}x+\begin{pmatrix}0\\
	v_{1}
	\end{pmatrix}y\right)=
	\end{equation}
	\begin{equation}
	=\underset{k=2}{\overset{N}{\sum}}c_{2}{\sqrt{N}}\varGamma_{k}^{2\delta_{1,\bullet}}\begin{pmatrix}\left(\dfrac{-2\lambda_{k}}{4\lambda_{k}^{2}+9}x^{2}+\dfrac{-6}{4\lambda_{k}^{2}+9}xy+\dfrac{2\lambda_{k}}{4\lambda_{k}^{2}+9}y^{2}+O\left(3\right)\right)V_{,k}\\
	\left(\dfrac{-\left(2\lambda_{k}^{2}+3\right)}{4\lambda_{k}^{2}+9}x^{2}+\dfrac{-2\lambda_{k}}{4\lambda_{k}^{2}+9}xy+\dfrac{-\left(2\lambda_{k}^{2}+6\right)}{4\lambda_{k}^{2}+9}y^{2}+O\left(3\right)\right)V_{,k}
	\end{pmatrix}\label{eq:Center_Manifold_2}
	\end{equation}

	If $f(x)=c_{3}x^{3}+O\left(x^{4}\right)$ then after scaling normalization we get:
	\begin{equation}
	\Psi\left(\begin{pmatrix}v_{1}\\
	0
	\end{pmatrix}x+\begin{pmatrix}0\\
	v_{1}
	\end{pmatrix}y\right)=
	\end{equation}
	\begin{equation}
	=\underset{k=2}{\overset{N}{\sum}}c_{3}{\sqrt{N}}^{2}\varGamma_{k}^{3\delta_{1,\bullet}}\begin{pmatrix}\left(\dfrac{-9\lambda_{k}^{2}-48}{9\lambda_{k}^{3}+64\lambda_{k}}x^{3}-\dfrac{18\lambda_{k}}{9\lambda_{k}^{3}+64\lambda_{k}}x^{2}y+\dfrac{-48}{9\lambda_{k}^{3}+64\lambda_{k}}xy^{2}-\dfrac{-6\lambda_{k}}{9\lambda_{k}^{3}+64\lambda_{k}}y^{3}+O\left(4\right)\right)V_{,k}\\
	\left(\dfrac{2\lambda_{k}}{9\lambda_{k}^{3}+64\lambda_{k}}x^{3}+\dfrac{-9\lambda_{k}^{2}-48}{9\lambda_{k}^{3}+64\lambda_{k}}x^{2}y+\dfrac{-6\lambda_{k}}{9\lambda_{k}^{3}+64\lambda_{k}}xy^{2}+\dfrac{-6\lambda_{k}^{2}-48}{9\lambda_{k}^{3}+64\lambda_{k}}y^{3}+O\left(4\right)\right)V_{,k}
	\end{pmatrix}\label{eq:Center_Manifold_3}
	\end{equation}
	
	Let us notice that for the matrix $M$ generated as in \eqref{eq:move_real_Lamda} we have that for $\Lambda \rightarrow \infty$ the 3rd order nonlinearity (odd) doesn't survive on the shape of the center manifold while the 2nd order nonlinearity (even) survives only on the second "layer" and only as even in both $x$ and $y$. This is not a general result but comes from the interplay between the vertical (module) and horizontal (M) connectivity. For instance similar calculations for coupled simple, unstructured $1D$ subunits gives a center manifold where even the even nonlinearities go to 0 when $\Lambda \rightarrow \infty$. 
	
	\section{Symmetries of $\alpha$'s and $\beta$'s}
	The above observation for the first order of approximation of $Psi$ is more general. In particular we can show that: $\alpha_{i}^{k}=\dfrac{P_{a,k,i}}{P}$, and $\beta_{i}^{k}=\dfrac{P_{\beta,k,i}}{P},$ where $P$ is a polynomial of $\lambda_{k}$ of order $k$ and $P_{\alpha,k,i}$ , $P_{\beta,k,i}$ polynomials of order $\leq k$. The inequality is streak except for $P_{\beta,k,i}$ when $k$ and $i$ are even. That would imply that for large $\Lambda$ only even nonlinearities survive by the operate $\Psi$ as even nonlinearities on the $y$ layer (the layer in which the original equations don't include nonlinearities). We will talk about this in details elsewhere but let us just notice that the coefficients are  just determined by determinants of submatrices of the pentadiagonal $A^{2}+\lambda_{k}A+id_{n}$ where:

	\noindent $\,$
	\noindent$A^{2}=$
	\begin{gather}\nonumber 	
	\left(\begin{array}{ccccccccccccc}
	-n & 0 & 1\cdot2 & 0\\
	0 & -n\cdot1-2\cdot\left(n-1\right) & 0 & 2\cdot3\\
	n\cdot\left(n-1\right) & 0 & -\left(n-1\right)\cdot2-3\cdot\left(n-2\right) & 0 & \ddots\\
	0 & \left(n-1\right)\cdot\left(n-2\right) & 0 & \ddots\\
	&  & \ddots &  & \ddots &  & \left(k-1\right)\cdot k\\
	&  &  &  &  & \ddots & 0 & k\cdot\left(k+1\right)\\
	&  &  &  & \left(n+2-k\right)\cdot\left(n+1-k\right) & 0 & -\left(n+1-k\right)\cdot k-\left(k+1\right)\cdot\left(n-k\right) & 0 & \left(k+1\right)\cdot\left(k+2\right)\\
	&  &  &  &  & \left(n+1-k\right)\cdot\left(n-k\right) & 0\\
	&  &  &  &  &  & \left(n-k\right)\cdot\left(n-k-1\right) &  &  &  & \ddots\\
	&  &  &  &  &  &  &  &  & \ddots & 0 & \left(n-2\right)\cdot\left(n-1\right) & 0\\
	&  &  &  &  &  &  &  & \ddots & 0 & -3\cdot\left(n-2\right)-\left(n-1\right)\cdot2 & 0 & \left(n-1\right)\cdot n\\
	&  &  &  &  &  &  &  &  & 3\cdot2 & 0 & -2\cdot\left(n-1\right)-n\cdot1 & 0\\
	&  &  &  &  &  &  &  & 0 & 2\cdot1 & 0 & -n
	\end{array}\right)
	1\end{gather}

	\section{Reduced Dynamics 3rd Order Approximation}
	
	\noindent Let $f\left(x\right)=c_{2}x^{2}+c_{3}x^{3}+O\left(x^{4}\right)$
	, then: 
	
	\noindent The reduced equations (\ref{eq:reduction}) are then up
	to order 3:
	
	\[
	\begin{pmatrix}v_{1}\\
	0
	\end{pmatrix}\dot{x}+\begin{pmatrix}0\\
	v_{1}
	\end{pmatrix}\dot{y}=W\begin{pmatrix}v_{1}\\
	0
	\end{pmatrix}x+W\begin{pmatrix}0\\
	v_{1}
	\end{pmatrix}y+
	\]
	
	\[
	+c_{2}\Pi_{c}\left(\begin{pmatrix}\left(V_{1,1}x+\underset{k=2}{\overset{N}{\sum}}V_{1,k}\left(\alpha_{0}^{k}x^{2}+\alpha_{1}^{k}xy+\alpha_{2}^{k}y^{2}\right)\right)^{2}\\
	\vdots\\
	\left(V_{N,1}x+\underset{k=2}{\overset{N}{\sum}}V_{N,k}\left(\beta_{0}^{k}x^{2}+\beta_{1}^{k}xy+\beta_{2}^{k}y^{2}\right)\right)^{2}\\
	0\\
	\vdots\\
	0
	\end{pmatrix}\right)+c_{3}\Pi_{c}\left(\begin{pmatrix}\left(V_{1,1}x\right)^{3}\\
	\vdots\\
	\left(V_{N,1}x\right)^{3}\\
	0\\
	\vdots\\
	0
	\end{pmatrix}\right)+O\left(4\right)
	\]

	So after scaling normalization we get the reduced dynamics equations:
	
	\begin{eqnarray*}
		\dot{x}&=&y+c_{2}{\sqrt{N}}\Gamma_{1}^{2\delta_{1,\bullet}}x^{2}+\left(c_{3}{\sqrt{N}}^{2}\varGamma_{1}^{3\delta_{1,\bullet}}+2c_{2}{\sqrt{N}}^{2}\underset{k=2}{\overset{N}{\sum}}\Gamma_{1}^{\delta_{1,\bullet}+\delta_{k,\bullet}}\alpha_{0}^{k}\right)x^{3}+\\
		&&+\left(2c_{2}{\sqrt{N}}^{2)}\underset{k=2}{\overset{N}{\sum}}\Gamma_{1}^{\delta_{1,\bullet}+\delta_{k,\bullet}}\alpha_{1}^{k}\right)x^{2}y++\left(2c_{2}{\sqrt{N}}^{2}\underset{k=2}{\overset{N}{\sum}}\Gamma_{1}^{\delta_{1,\bullet}+\delta_{k,\bullet}}\alpha_{2}^{k}\right)xy^{2}+O\left(4\right)\\
		\dot{y}&=&-x
	\end{eqnarray*}
	
	\noindent \underline{Proposition 3: }Let $M_{0}$ with $\lambda_{1}=0,$
	$Re\left(\lambda_{i}\right)<0$ for $i>1$ and $M\left(\mu\right)=M_{0}+\mu V_{,1}W_{1,}$
	. The system (\ref{eq:leyerx}) undergoes an Andronov-Hopf bifurcation at $\mu=0$ if $a=c_{3}\dfrac{3}{8}{\sqrt{N}}^{2}\varGamma_{1}^{3\delta_{1,\bullet}}-c_{2}^{2}{\sqrt{N}}^{2}\underset{k=2}{\overset{N}{\sum}}\Gamma_{1}^{\delta_{1,\bullet}+\delta_{k,\bullet}}\varGamma_{k}^{2\delta_{1,\bullet}}\dfrac{\lambda_{k}}{4\lambda_{k}^{2}+9}\neq0.$
	
	The bifurcation is subcritical (supercritical) iff $a>0$ ($a<0$).
	
	\noindent \underline{Proof sketch:} Proposition 2b $\Rightarrow$
	nonhyperbolicity and transversality condition. For the system:
	\begin{eqnarray*}
		\dot{x}&=&y+\tilde{f}\left(x,y\right)\\
		\dot{y}&=&-x+\tilde{g}\left(x,y\right)
	\end{eqnarray*}
	
	\noindent the genericity condition parameter is :
	
	\noindent $a=\dfrac{1}{16}\left(\tilde{f}_{xxx}+\tilde{f}_{xyy}+\tilde{g}_{xxy}+\tilde{g}_{yyy}\right)+\dfrac{1}{16\omega}\left(\tilde{f}_{xy}\left(\tilde{f}_{xx}+\tilde{f}_{yy}\right)-\tilde{g}_{xy}\left(\tilde{g}_{xx}+\tilde{g}_{yy}\right)-\tilde{f}_{xx}\tilde{g}_{xx}+\tilde{f}_{yy}\tilde{g}_{yy}\right)=\dfrac{1}{16}\left(\tilde{f}_{xxx}+\tilde{f}_{xyy}\right)=c_{3}\dfrac{3}{8}{\sqrt{N}}^{2}\varGamma_{1}^{3\delta_{1,\bullet}}-c_{2}^{2}{\sqrt{N}}^{2}\underset{k=2}{\overset{N}{\sum}}\Gamma_{1}^{\delta_{1,\bullet}+\delta_{k,\bullet}}\varGamma_{k}^{2\delta_{1,\bullet}}\dfrac{\lambda_{k}}{4\lambda_{k}^{2}+9}\blacksquare$
	
	In particular let us notice that this operator can drastically change the whole dynamics and create a subcritical bifurcation from subunits of supercritical bifurcation.

	Notice: if $M$ is normal then  $-\sqrt{N}^{\left(3-1\right)}c_{2}^{2}\underset{k=2}{\overset{N}{\sum}}\Gamma_{1}^{\delta_{1,\bullet}+\delta_{k,\bullet}}\varGamma_{k}^{2\delta_{1,\bullet}}\dfrac{\lambda_{k}}{4\lambda_{k}^{2}+9}>0$.
	$\xrightarrow{N\rightarrow\infty}3c_{2}^{2}\dfrac{\lambda_{k}}{4\lambda_{k}^{2}+9}$,
	and the naive approximation of defining $\Phi$ should be appropriate
	for $c_{2}\ll\sqrt{\Lambda}$. 
	
\begin{figure}[]
	\begin{centering}
		\includegraphics[scale=0.7]{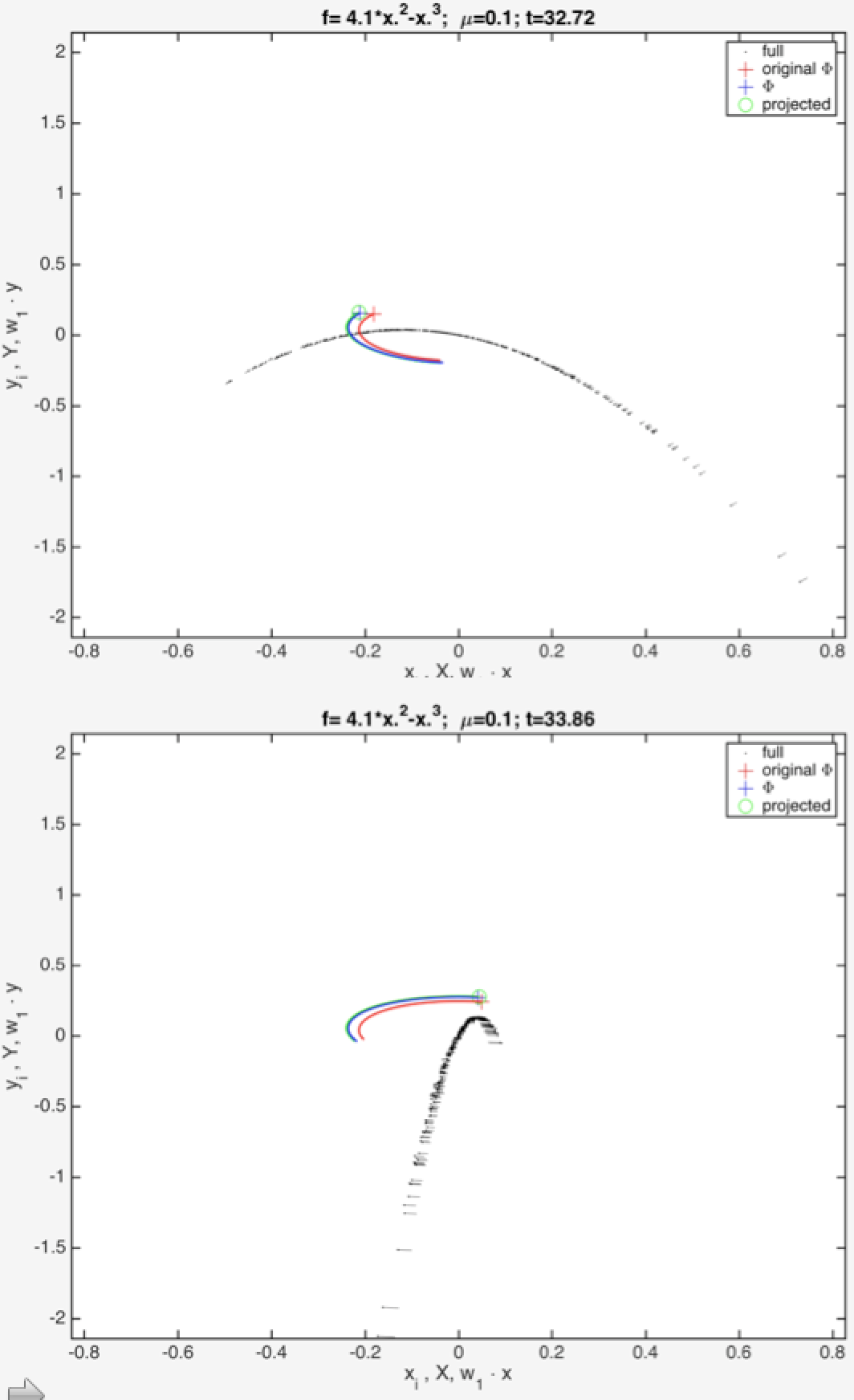}
		\par\end{centering}
	
	\begin{centering}
		\caption{\label{fig:phi_correction} \textbf{Corrected $\Phi$} As in Fig. \ref{fig:naive_phi_fail} adding in blue color the correction using the second order terms in the approximation of $\Phi$ to capture the reduced dynamics.}
		\par\end{centering}
		
	\end{figure}

	Now we can use the nonlinearities of the newly defined reduced dynamics  to define our new $\Phi$. Which is shown as a blue cross in movie \url{https://drive.google.com/open?id=1NN0wju5Dl_VBD0JkL-PBHvfmJkWKdh8F}. Snapshots of the simulation movie are shown in (Fig. \ref{fig:phi_correction}).
	
	 \section{Numerical Ansatz for Approximating $\Psi$ When $f\left(x\right)=c_{2}x^{2}+c_{3}x^{3}$ and $M$ is Normal and $\lambda_{k}\approx-d$ for $k\geq2$.}
	
	Let's use the 2nd order nonlinearities of the center manifold (at
	the bifurcation point) as a guess of the residuals of individual units for small $\lambda_{1}>0$ . 
	
	\[
	x_{i}={\sqrt{N}}V_{i,1}X_{\lambda_{1}}+\underset{k=2}{\overset{N}{\sum}}{\sqrt{N}}V_{i,k}X_{\lambda_{k}}\approx {\sqrt{N}}V_{i,1}X_{\lambda_{1}}+\underset{k=2}{\overset{N}{\sum}}{\sqrt{N}}V_{i,k}c_{2}{\sqrt{N}}\varGamma_{k}^{2\delta_{1,\bullet}}\left(\dfrac{-2\lambda_{k}}{4\lambda_{k}^{2}+9}X_{\lambda_{1}}^{2}+\dfrac{-6}{4\lambda_{k}^{2}+9}X_{\lambda_{1}}Y_{\lambda_{1}}+\dfrac{2\lambda_{k}}{4\lambda_{k}^{2}+9}Y_{\lambda_{1}}^{2}\right)
	\]
	
	\noindent If for $k\geq2$ ,$\lambda_{k}\approx-d$ then:
	
	\[
	x_{i}\approx {\sqrt{N}}V_{i,1}X_{\lambda_{1}}+c_{2}\underset{k=2}{\overset{N}{\sum}}{\sqrt{N}}^{2}V_{i,k}\varGamma_{k}^{2\delta_{1,\bullet}}\left(\dfrac{2d}{4d^{2}+9}X_{\lambda_{1}}^{2}+\dfrac{-6}{4d^{2}+9}X_{\lambda_{1}}Y_{\lambda_{1}}+\dfrac{-2d}{4d^{2}+9}Y_{\lambda_{1}}^{2}\right) 
	\]
	
	\noindent If $M$ is normal then : $\underset{k=2}{\overset{N}{\sum}}V_{i,k}\varGamma_{k}^{2\delta_{1,\bullet}}=\underset{k=2}{\overset{N}{\sum}}V_{i,k}\underset{\varphi=1}{\overset{N}{\sum}}V_{\varphi,k}V_{\varphi,1}^{2}=\underset{\varphi=1}{\overset{N}{\sum}}V_{\varphi,1}^{2}\left(\delta_{i,\varphi}-V_{i,1}V_{\varphi,1}\right)=V_{i,1}^{2}-V_{i,1}\Gamma_{1}^{2\delta_{1,\bullet}}$.

	\noindent Since $\left(\sqrt{N}\right)^{2}V_{i,1}\Gamma_{1}^{2\delta_{1,\bullet}}\xrightarrow{N\rightarrow\infty}0$.
	We get:
	\begin{eqnarray*}
		x_{i}&\approx& \left(\sqrt{N}\right)V_{i,1}X_{\lambda_{1}}+c_{2}\left(\sqrt{N}\right)^{2}V_{i,1}^{2}\left(\frac{2d}{4d^{2}+9}X_{\lambda_{1}}^{2}+\frac{-6}{4d^{2}+9}X_{\lambda_{1}}Y_{\lambda_{1}}+\frac{-2d}{4d^{2}+9}Y_{\lambda_{1}}^{2}\right)\\
		y_{i}&\approx& V_{i,1}\left(\sqrt{N}\right)Y_{\lambda_{1}}+c_{2}\left(\sqrt{N}\right)^{2}V_{i,1}^{2}\left(\frac{-\left(2d^{2}+3\right)}{4d^{2}+9}X_{\lambda_{1}}^{2}+\frac{2d}{4d^{2}+9}X_{\lambda_{1}}Y_{\lambda_{1}}+\frac{-\left(2d^{2}+6\right)}{4d^{2}+9}Y_{\lambda_{1}}^{2}\right)
	\end{eqnarray*}
	
\begin{figure}[]
		\begin{centering}
			\includegraphics[scale=0.7]{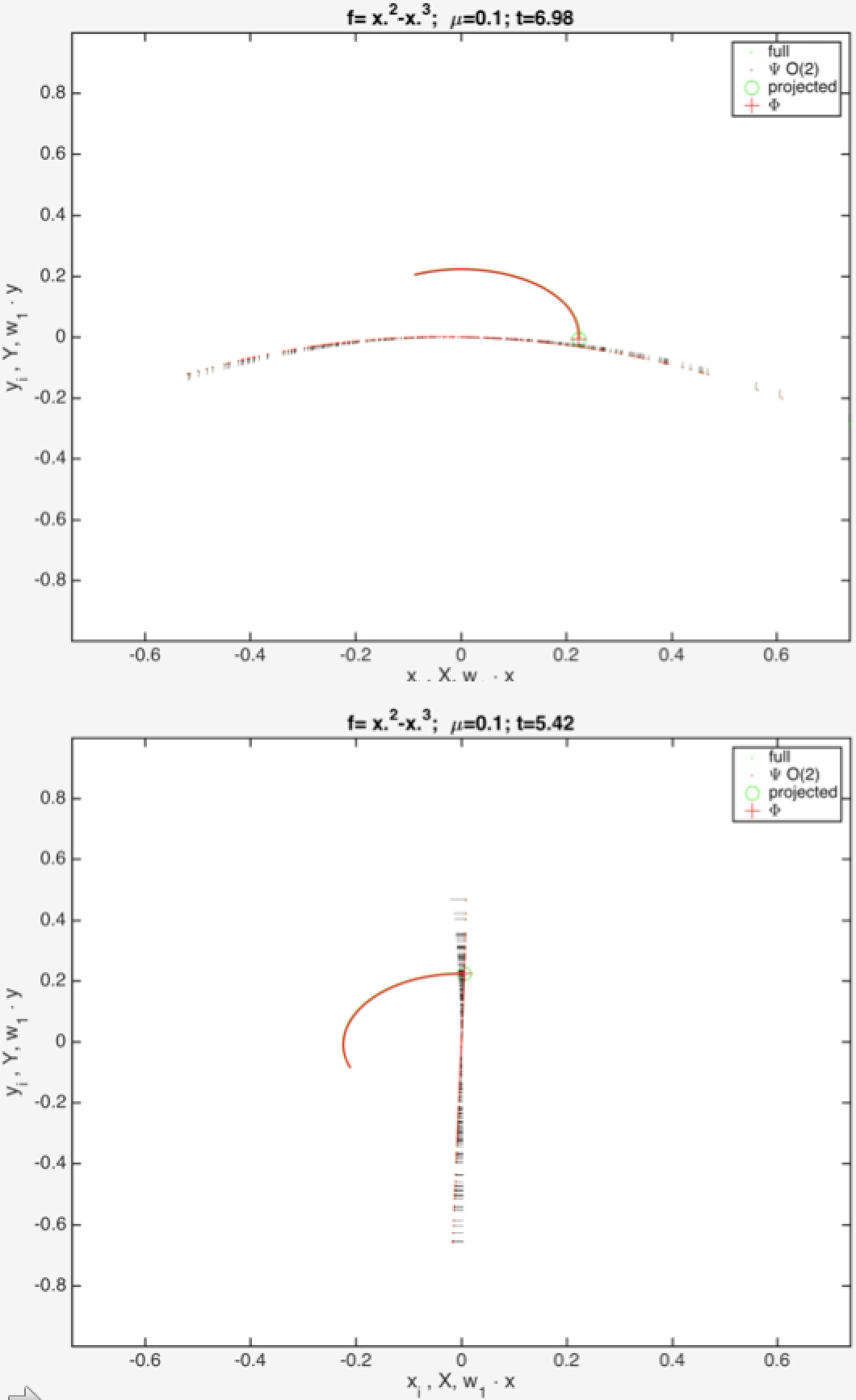}
			\par\end{centering}
		
		\begin{centering}
			\caption{\label{fig:psi_O(2)} \textbf{Corrected $\Psi$ with order 2 terms} As in Fig. \ref{fig:naive2} adding in red color the second order terms in the of approximation $\Psi$ to capture the dynamics of individual units from the reduced dynamics. }
			\par\end{centering}
		
	\end{figure}

	The corrected $\Psi$ is shown in green color in the movie: 
	\url{https://drive.google.com/open?id=1GAUPMhUs04eMjNSbMxjPw8MK2kIGlTZG} . Simulation is in red color. Snapshots of the simulation movie are shown in (Fig. \ref{fig:psi_O(2)}).
	
\begin{figure}[]
		\begin{centering}
			\includegraphics[scale=0.7]{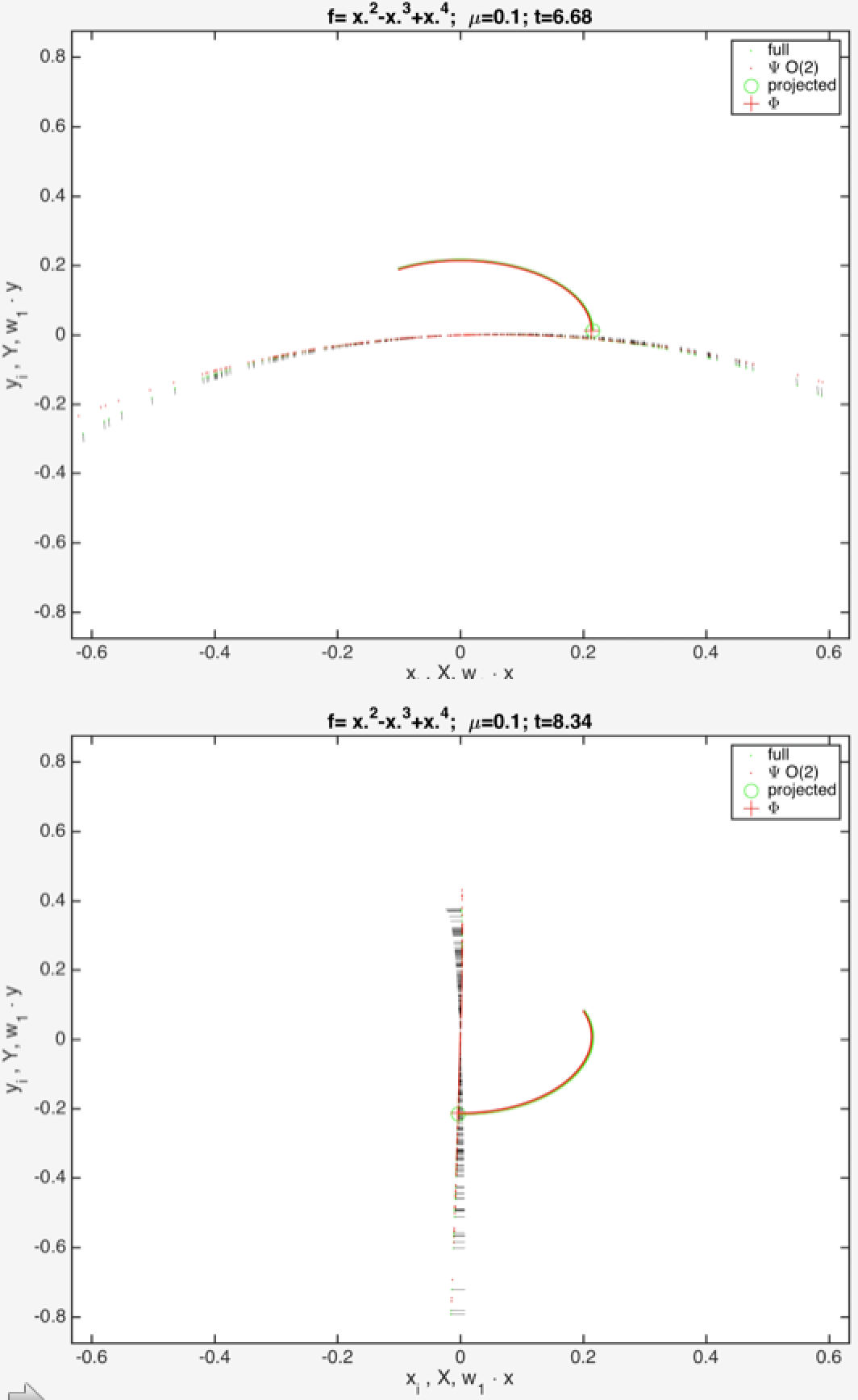}
			\par\end{centering}
		
		\begin{centering}
			\caption{\label{fig:psi_O(2)_fail} \textbf{Corrected $\Psi$ with Order 2 Terms Fails to Capture 4th Order Terms} As in (Fig. \ref{fig:psi_O(2)}) with $f(x)=x^2-x^3+x^4$. Only second order terms of $\Psi$ are not enough. }
			\par\end{centering}
		
	\end{figure}

	We see that higher order even nonlinearities also survive on $\Psi$ . \url{https://drive.google.com/open?id=1mb3QK2JGSVxShdKyUA_0bNDhOgSAk_tA}. Snapshots of the simulation movie are shown in (Fig. \ref{fig:psi_O(2)_fail}).

	\section{Numerical Ansatz for Approximating $\Psi$ When $f\left(x\right)=c_{2}x^{2}+c_{3}x^{3}+c_{4}x^{4}+c_{5}x^{5}$ and $M$ is Normal and $\lambda_{k}\approx-d$ for $k\geq2$.}
	
	Assuming that $O(2)*O(2)$ in (\ref{eq:approx}) disappear in the
	limit and using the highly simplified guess from (\ref{eq:Center_Manifold_Parameters})
	and adding linear addition $O(2)+O(4)$ we get the approximations on
	the last movie. In the particular case described here this can intuitively be seen to be an appropriate simplification for $\Lambda \rightarrow \infty$ due to the nonlinearities only surviving on the x layer on $\Phi$ and on the y layer on $\Psi$.
	
\begin{figure}[]
	\begin{centering}
		\includegraphics[scale=0.7]{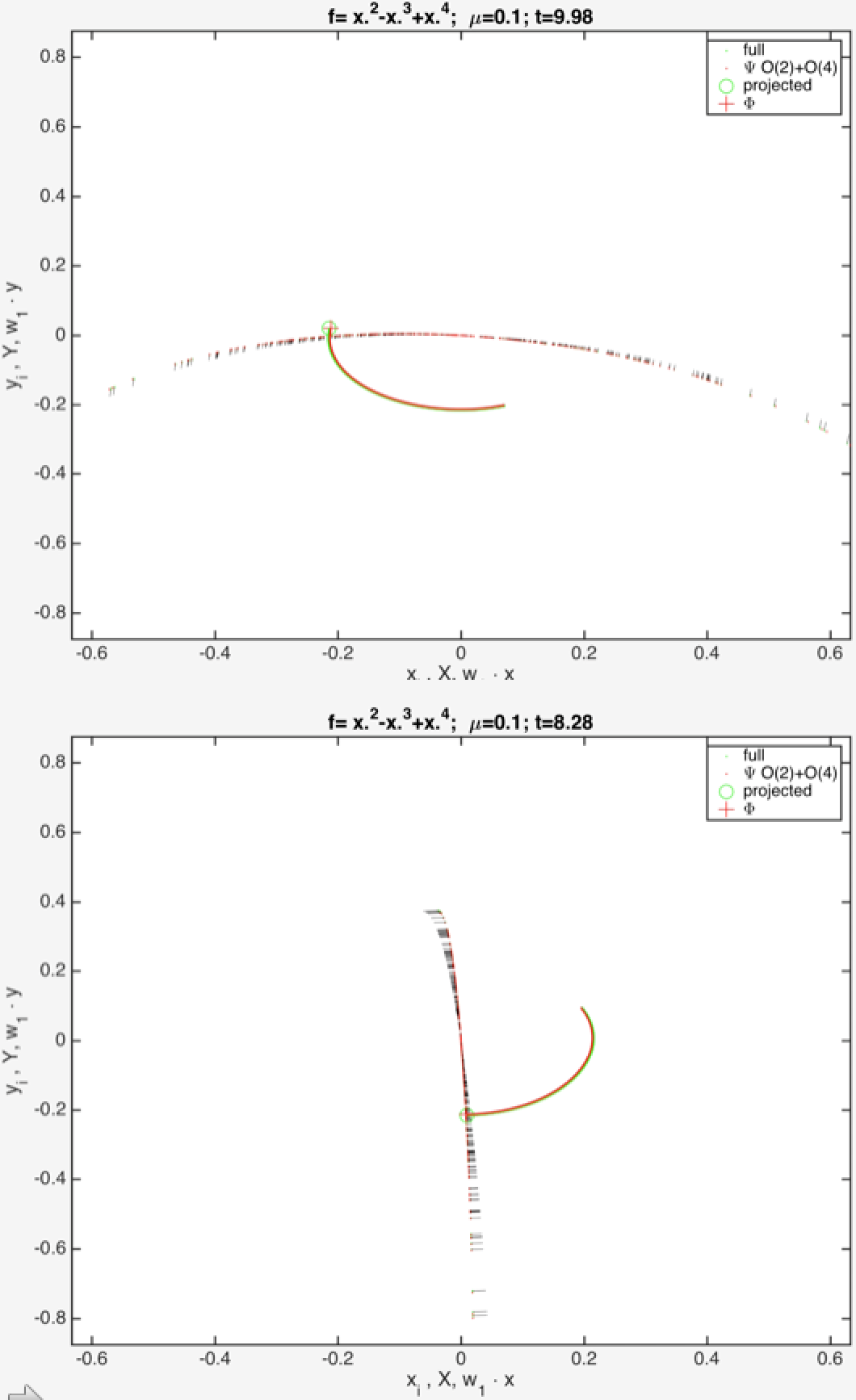}
		\par\end{centering}
	
		\begin{centering}
		\caption{\label{fig:psi_symplified_O(2)+O(4)} \textbf{Simplified $\Psi$ Containing Order 2 and 4 Terms} As in Fig. \ref{fig:psi_O(2)_fail} but adding terms of order 4 in the approximation of $\Psi$. }
		\par\end{centering}
	
	\end{figure}

	The corrected $\Psi$ is shown in green color in the movie:  \url{https://drive.google.com/open?id=1pXFWluButTHG1jReFUm6Rv9J1ihzrNKe}.
	Snapshots of the simulation movie are shown in (Fig. \ref{fig:psi_symplified_O(2)+O(4)}).
	
	\section{Complex-Conjugate Pair of Eigenvalues, Normal Matrix}
	
	As shown in (Fig. \ref{fig:eigenvalue_diagram}), the second-simplest case is to control two complex-conjugate eigenvalues, as they can be controlled by a single parameter, their real part, while leaving their imaginary part unchanged. 
	
	We follow the previous section in using $M_a$, a normal matrix, as our base connectivity structure. Now we choose a pair of complex-conjugate eigenvalues. The corresponding eigenvectors are also complex conjugate. Choosing as ${ e}$ the one corresponding to the eigenvalue with positive imaginary part, the 2D subspace is spanned both by ${ e}$ and $ \bar{{ e}}$, or by ${ e}_R \equiv \Re{ e}$ and ${ e}_I \equiv \Im{\ e}$. Corollary 2 implies that the whole system bifurcates with 2 pairs of complex conjugate eigenvalues crossing the imaginary axis with frequencies  $\dfrac{\sqrt{\omega^{2}+4} \pm \omega}{2}$ .
	
	Following the previous section, and remembering that now the left eigenvector is the complex conjugate of the right eigenvector, ${ e}^\dag = \bar{{ e}}^t $ we define the coarse-grained variable $X$, now complex, as the normalized projection on ${ e}$:
	$$X\sqrt{\dfrac{N}{2}} = \left({ e}_{R}-{ e}_{I}i\right)\cdot { x}$$
	
	projection
	on critical mode, $X=X_{R}+X_{I}i$

	Then the linear order approximation of the invariant manifold for the system:
	\begin{eqnarray*}
		\dot{x}_{i} & = & y_{i}-x_{i}^{\alpha}+\underset{j}{\sum}M_{ij}x_{j}\\
		\dot{y}_{i} & = & -x_{i}
	\end{eqnarray*}
	gives us the linear order of $\Psi$ 
	\begin{eqnarray*}
		{ x} & \approx & \sqrt{\dfrac{N}{2}}\left(X\left({ e}_{R}+ i { e}_{I} \right)+\bar{X}\left({ e}_{R}-i { e}_{I}\right)\right)
	\end{eqnarray*}
	
	\noindent The naive approximation of $\Psi$ is given by :
	\begin{eqnarray*}
	\end{eqnarray*}

	\begin{eqnarray*}
		x & = & \sqrt{\dfrac{N}{2}}\left(\left(X_{R}+X_{I}i\right)\left(e_{R}+e_{I}i\right)+\left(X_{R}-X_{I}i\right)\left(e_{R}-e_{I}i\right)\right)=\sqrt{2N}\left(X_{R}e_{R}-X_{I}e_{I}\right)\\
		y & = & \sqrt{\dfrac{N}{2}}\left(\left(Y_{R}+Y_{I}i\right)\left(e_{R}+e_{I}i\right)+\left(Y_{R}-Y_{I}i\right)\left(e_{R}-e_{I}i\right)\right)=\sqrt{2N}\left(Y_{R}e_{R}-Y_{I}e_{I}\right)
	\end{eqnarray*}
	
	\noindent Then taking dot product with $\sqrt{\dfrac{N}{2}}e_{R}$ and $-\sqrt{\dfrac{N}{2}}e_{L}$,
	and ignoring projections other than in $e_{R}$ and $e_{I}$ the system:
	
	\begin{eqnarray*}
		\dot{x}_{i} & = & y_{i}-x_{i}^{\alpha}+\underset{j}{\sum}M_{ij}x_{j}\\
		\dot{y}_{i} & = & -x_{i}
	\end{eqnarray*}

	$\,$
	
	\noindent gives us :
	
	\begin{eqnarray*}
		\dot{X}_{R} & = & Y_{R}+\lambda_{R}X_{R}-\lambda_{I}X_{I}-\sum\sqrt{\dfrac{2}{N}}e_{R,i}\sqrt{2N}^{\alpha}\left(X_{R}e_{R,i}-X_{I,i}e_{I,i}\right)^{\alpha}\\
		\dot{Y}_{R} & = & -X_{R}\\
		\dot{X}_{I} & = & Y_{I}+\lambda_{I}X_{R}+\lambda_{R}X_{I}-\sum\sqrt{\dfrac{2}{N}}\left(-e_{I,i}\right)\sqrt{2N}^{\alpha}\left(X_{R}e_{R,i}-X_{I,i}e_{I,i}\right)^{\alpha}\\
		\dot{Y}_{I} & = & -X_{I}
	\end{eqnarray*}

	\noindent for $\alpha$ odd:

	\noindent$\sum\sqrt{\dfrac{2}{N}}e_{R,i}\sqrt{2N}^{\alpha}\left(X_{R}e_{R,i}-X_{I,i}e_{I,i}\right)^{\alpha}=\underset{\varphi =0}{\overset{\alpha}{\sum}}\begin{pmatrix}\alpha\\
	\varphi
	\end{pmatrix}\dfrac{1}{N}\sqrt{2N}^{\alpha+1}\underset{i}{\sum}e_{R,i}^{\alpha-\varphi}\left(-e_{I,i}\right)^{\varphi}=$
	
	$X_{R}\underset{\varphi =0}{\overset{\frac{\alpha-1}{2}}{\sum}}\begin{pmatrix}\alpha\\
	2\phi
	\end{pmatrix}\left(\alpha-2\phi+1\right)!!\left(2\phi\right)!!X_{R}^{\alpha-2\phi-1}X_{I}^{2\phi}=\alpha!!X_{R}\left(X_{R}^{2}+X_{I}^{2}\right)^{\frac{\alpha-1}{2}}$
	
	\noindent since $\begin{pmatrix}\alpha\\
	2\phi
	\end{pmatrix}\left(\alpha-2\phi+1\right)!!\left(2\phi\right)!!=\alpha!!\begin{pmatrix}\frac{\alpha-1}{2}\\
	\phi
	\end{pmatrix}$

	\noindent similarly $\sum\sqrt{\dfrac{2}{N}}\left(-e_{I,i}\right)\sqrt{2N}^{\alpha}\left(X_{R}e_{R,i}-X_{I,i}e_{I,i}\right)^{\alpha}=\alpha!!X_{I}\left(X_{R}^{2}+X_{I}^{2}\right)^{\frac{\alpha-1}{2}}$
	
	$\,$
	
	\noindent So:
	
	\begin{eqnarray*}
		\dot{X} & = & Y+\lambda X-\alpha!!\left\Vert X\right\Vert ^{\alpha-1}X\\
		\dot{Y} & = & -X
	\end{eqnarray*}

	For $\alpha = 3$ this is a complex generalized van der Pol equation. 
	\begin{eqnarray*}
		\dot{X} & = &  (\lambda - |X|^2) X + Y \\
		\dot{Y} & = & -X
	\end{eqnarray*}
	Notice the $X$ equation is now a Hopf bifurcation normal form in a single complex variable, where ${ Imag}(\lambda)$ is the rotational frequency in the $X$ complex plane, and ${ Real}(\lambda)$ is the control parameter, while the whole $XY$ plane equation uses the same ${ Real}(\lambda)$ as a control parameter for a Hopf bifurcation in the real $XY$ plane. Thus this equation contains 2 distinct frequencies, $\omega = 1$ from the base equation, and $\omega = { Imag}(\lambda)$, the frequency of the eigenvector, from the complex $X$ equation, and both cross the threshold of stability at the same time. Therefore, in this collective mode, a codimension-2 transition becomes generic, and the system in principle can bifurcate directly from a fixed point onto a torus. In this sense it is similar to a Takens-Bogdanoff normal form, only that it has two distinct and independent frequencies. 
	
	Thus the reduced equations are non-generic with respect to perturbations of the individual units in the full system in two ways: first, all even terms are destroyed, and second it can undergo a direct transition from a fixed point to a torus. 

	\section{Conclusion}
	
	We have proposed a statistical center manifold reduction technique based on the three foundational center manifold theorems: existence, stability, and approximation \cite{Kelley1967,Carr2012}. Utilizing these theorems we have constructed two key operators. The first operator $\Phi$ connects the statistics of the vector field in the full-dimensional space to the statistics of the reduced vector field in the low-dimensional space. The second operator $\Psi$ lifts the statistical properties of the trajectories produced by these reduced vector fields to the trajectories evolved by the original vector fields. For any given system, these two operators are sufficient to verify that our center manifold dimensionality reduction approach is well defined (the diagram in Fig. \ref{fig:schematics2} commutes).

	In this paper, we have considered a specific example of a simple network of randomly connected, spontaneous oscillators. We have shown that our center manifold dimensionality reduction approach is well defined, successfully capturing the statistical properties of the individual units in the full system. The operators $\Phi$ and $\Psi$ depend crucially on the structure of the eigenvectors (through the random variables $\Gamma$'s). While we now know quite a lot about the statistical nature of the eigenvalues of random matrices \cite{Tao2008circular}, little is known about the statistical structure of the eigenvectors of these matrices \cite{Chalker1998Eigenvector, Tao2012Random, O'Rourke2016Eigenvectors}. In the considered example, the reduction operator $\Phi$ coarse grains out the even nonlinearities present in the full system. But even though the even nonlinearities are coarse-grained out, we have found that they can still leak to higher order odd nonlinearities, which can have drastic effects on the dynamics as a whole, and for example, create a subcritical Hopf bifurcation from subunits of supercritical Hopf bifurcations. Moreover, using simplifying assumptions in the approximation theorem, we were able to construct a lift operator $\Psi$ for higher order nonlinearities. This lift operator has interesting properties like in the limit of $d\rightarrow\inf$, the even nonlinearities survive only in the y layer and the odd nonlinearities vanish entirely.
	
	We have explored several other systems and it appears that the property of even nonlinearities being coarse-grained out in the reduced dynamics generally holds for a class of similar systems although the exact nature of this class is not yet clear. Furthermore, it should be easily to generalize to other structures such as sparsity by calculating new values for the $\Gamma$'s. Both these subjects should be a target of future research.

\end{document}